\documentclass{amsart}

\usepackage{amsmath,amsfonts,amssymb,amsthm,amscd,latexsym,euscript,bbm}

   \newtheorem{lemma}{Lemma}[section]  
   \newtheorem{kor}[lemma]{Corollary}        
   \newtheorem{satz}[lemma]{Theorem}         

\theoremstyle{definition}
\newtheorem{Def}[lemma]{Definition}                                          

    \newtheorem{Bem}[lemma]{Remark}        
    \newtheorem{Beisp}[lemma]{Example}     
    \newtheorem{Ann}[lemma]{Assumptions}   

\usepackage[all]{xy}

\newenvironment{beweis}{\noindent\textbf{Proof:}}{\par \hfill $\Box$ \hspace{1cm} \par}%

\newenvironment{proofof}[1]{\renewcommand\proofname{#1} {\noindent\bf{Proof of #1:}}}%
{\par  \hfill {\bf{$\Box $}} \hspace{1cm}       \par}

\newcounter{zahl}%
\newenvironment{punkt}{\begin{list}{{\rm{(\roman{zahl})}}}%
    {\usecounter{zahl}%
     \setlength{\leftmargin}{0pt} \setlength{\itemindent}{4pt} \setlength{\topsep}{2pt} \setlength{\parsep}{2pt} }}%
    {\end{list}}%

\newcommand{\co}{\colon\thinspace}
\newcommand{\mc}[1]{\ensuremath{\mathcal{#1}}}
\newcommand{\sm}[2]{\ensuremath{#1\smash[b]{\wedge\,}#2}}
\newcommand{\aus}{\raisebox{1pt}{\ensuremath{\,{\scriptstyle\in}\,}}}%

\newcommand{\toh}[1]{\ensuremath{\stackrel{#1}{\rightarrow}}}%
\newcommand{\colim}{\operatorname*{colim}}
\newcommand{\sk}{\operatorname{sk}}
\newcommand{\cosk}{\operatorname{cosk}}
\newcommand{\holim}{\operatorname*{holim\,}}
\newcommand{\hocolim}{\operatorname*{hocolim\,}}
\newcommand{\Rlim}{\operatorname*{lim^1}}
\newcommand{\frei}{\,\_\!\_\,}
\newcommand{\id}{\operatorname{id}}
\newcommand{\Tot}{\operatorname{Tot}}
\newcommand{\inj}{\ensuremath{\hookrightarrow}}
\newcommand{\Fac}{\operatorname{Fac}}
\newcommand{\Real}{\operatorname{Real}}

\newcommand{\kl}{\ensuremath{\raisebox{1pt}{\{}}}
\newcommand{\kr}{\ensuremath{\raisebox{1pt}{\}}}}
\newcommand{\hel}{\ensuremath{\raisebox{1pt}{[}}}
\newcommand{\her}{\ensuremath{\raisebox{1pt}{]}}}
\newcommand{\hrl}{\ensuremath{\raisebox{1pt}{(}}}
\newcommand{\hrr}{\ensuremath{\raisebox{1pt}{)}}}
\newcommand{\ho}[1]{\ensuremath{H\!o({#1})}}
\newcommand{\bth}{\raisebox{1pt}{\ensuremath{\,\scriptstyle \geq\,}}}%
\newcommand{\sth}{\raisebox{1pt}{\ensuremath{\,\scriptstyle \leq\,}}}%
\newcommand{\uber}{\raisebox{1pt}{\ensuremath{\,{\scriptstyle > }\,}}}%
\newcommand{\unter}{\raisebox{1pt}{\ensuremath{\,{\scriptstyle < }\,}}}%
\newcommand{\winkel}{\ulcorner\hspace{-3pt}\raisebox{1pt}{\ensuremath{\cdot}}}%
\newcommand{\pushout}{\raisebox{-3pt}{\ensuremath{\displaystyle \winkel}}}%
\newcommand{\ol}[1]{\overline{#1}}
\newcommand{\wt}[1]{\widetilde{#1}}
\newcommand{\stabhom}[2]{\ensuremath{\hel#1,#2\her}}%
\newcommand{\Hom}[3]{\ensuremath{{\rm Hom}_{#1}\hrl#2,#3\hrr}}%
\newcommand{\Ext}[4]{\ensuremath{{\rm Ext}_{#1}^{#2}\hrl#3,#4\hrr}}%

\newcommand{\naturalpi}[3]{\ensuremath{\pi_{#1}^{\natural}\hrl #2,#3\hrr}}%

\newcommand{\map}{\ensuremath{{\rm map}}}
\renewcommand{\hom}{\ensuremath{{\rm hom}}}

\newcommand{\diag}[2]{ \begin{align} \begin{split} \xymatrix{#1} \end{split} \label{#2} \end{align}}%

\newcommand{\diagr}[1]{ \begin{equation*} \xymatrix{#1} \end{equation*}}%

\numberwithin{equation}{section}        
\newcommand{\Ref}[1]{\hrl\ref{#1}\hrr}        

\begin{document}

\SelectTips{cm}{10}

\title  [Interpolation categories]
        {Interpolation categories for homology theories}
\author{Georg Biedermann}

\address{Department of Mathematics, Middlesex College, The University of Western Ontario, London, Ontario N6A 5B7, Canada}

\email{gbiederm@uwo.ca}



\date{\today}
\dedicatory{F\"ur Ralf F\"utterer}
\commby{}

\begin{abstract}
For a homological functor from a triangulated category to an abelian category satisfying some technical assumptions we construct a tower of interpolation categories. These are categories over which the functor factorizes and which capture more and more information according to the injective dimension of the images of the functor. The categories are obtained by using truncated versions of resolution model structures. Examples of functors fitting in our framework are given by every generalized homology theory represented by a ring spectrum satisfying the Adams-Atiyah condition. The constructions are closely related to the modified Adams spectral sequence and give a very conceptual approach to the associated moduli problem and obstruction theory. As application we establish an isomorphism between certain E(n)-local Picard groups and some Ext-groups. 
\end{abstract}

\maketitle

\section{Introduction} 
\label{section:intro}

Algebraic topology, or more precisely homotopy theory, is the study of geometric objects up to weak homotopy equivalence by translating the geometrical or homotopical information into algebraic data.  
The mathematical device to do this are functors from a topological category to an algebraic category.
Computationally very useful are homological functors or homology theories, which are functors satisfying excision. 
Example in classical algebraic topology are abundant: singular homology, $K$-theory and cobordism, and more. For the purpose of this article we adopt the general definition, that a homological functor $F\co\mc{T}\to\mc{A}$ from a triangulated category \mc{T} to an abelian category \mc{A} is an additive functor, which maps distinguished triangles to long exact sequences. Here \mc{A} is graded in the following sense: there is a self equivalence $\hel 1\her$ of \mc{A} called shift, such that
    $$ F\hrl\Sigma X\hrr\cong\hrl FX\hrr\hel 1\her $$
via a natural isomorphism. The functors we consider have to satisfy some standard technical conditions \ref{Annahmen für F} met by most examples in topology.

For such a functor we would like to study questions of the following kind: Given an object $A$ in the abelian target category \mc{A}, does there exist an object $X$ in the triangulated source category \mc{T} together with an isomorphism $FX\cong A$ and if yes, how many different objects exist? We can ask the same question for morphisms.
This is the realization or moduli problem for the homological functor $F$.  

This article contributes a theory of interpolation categories for $F$ and an obstruction calculus for lifting objects or morphisms.
The categories are intended to give a very conceptual approach to realizations and moduli problems of this functor and to the associated obstruction theory. They interpolate in a precise sense between the topological source and the algebraic target category. The obstruction groups will lie in the $E_2$-term of the $F$-based Adams spectral sequence.

The way, we set up the obstruction calculus, follows the philosophy of \cite{BlDG:pi-algebra} and \cite{GoHop:moduli}. We use injective resolutions in \mc{A} and try to realize them step by step as cosimplicial objects over \mc{M}, where \mc{M} is a stable model category having \mc{T} as its homotopy category. As we go along by gluing certain layers \ref{L(N,n)} to our potential liftings to kill their higher cohomology, we run into difficulties giving us the obstruction groups, see theorems \ref{Ext^(n+2,n) und Objekte} and \ref{Ext^(n+1,n) und Objekte}. They can be described as certain mapping spaces \ref{Kohomologie} out of these layers. 

In order to be able to identify cosimplicial objects stemming from different resolutions, we consider resolution model structures on the category $c\mc{M}$ of cosimplicial objects. They were invented in \cite{DKSt:E2} and used for exactly this purpose.
Recently Bousfield in \cite{Bou:cos} has given a very general and elegant treatment of resolution model structures, which exhibits all previous instances as special cases.

In \cite{Biedermann:truncated} I truncate these resolution model structures and this article is based on the results there. Weak equivalences are now given by maps, which induce isomorphisms on cohomology just up to degree $n$, \ref{Def. n-G-Struktur}. The successive stages in our realization process are now given as cofibrant objects in these truncated resolution model structures. For an object $X$ in \mc{M} they are given by the skeletons of a Reedy cofibrant replacement of the constant cosimplicial object over $X$ and fit into a filtration, called the Postnikov cotower \ref{Ko-Postnikov-Turm}. This organizes the obstruction calculus nicely and enables us to define interpolation categories \ref{Interpolationskategorien}.
We carry out in full detail the obstruction calculus for realizing maps, theorems \ref{Ext^(n+1,n) und Morphismen} and \ref{Ext^(n,n) und Morphismen}, and prove in \ref{Baues-Turm}, that there is a full-fleshed theory of interpolation categories as axiomatized in \cite{Baues:combinatorial}, which keeps track of how obstruction behave under composition. 
We also obtain a description \ref{Adams-Differential} of the Adams differential $d_n\co E_n^{0,0}\hrl X,Y\hrr\to E_n^{n,n-1}\hrl X,Y\hrr$ in terms of our co-k-invariants.

When we finally have realized an object as an $\infty$-stage, we apply $\Tot\co c\mc{M}\to\mc{M}$. If the relevant spectral sequence converges, which happens in exactly those cases, when $F$-localization and $F$-completion coincide, we obtain an object in \mc{M}, which is a realization of the object in \mc{A}, we started with, see \ref{möglichweise isomorph} and \ref{die Abbildung Real_8(A) -> Real(A)}.

As in \cite{BlDG:pi-algebra} and \cite{GoHop:moduli} we study also moduli spaces of objects and morphisms, but with the truncated model structures at hand the proofs become a lot easier and shorter.

Finally we apply our obstruction calculus to the problem of determining certain $E\hrl n\hrr$-local Picard groups. In theorem \ref{Pic=Ext} we establish an isomorphism between them and some Ext-groups for a certain range of $n$ and $p$ building on and extending results from \cite{Hov-Sad:invertible}. However this result just uses the obstruction calculus and not the full result on the existence of interpolation categories. \\

\section{Resolution model structures} 
\label{section:E2-Strukturen}
Let \mc{M} be a model category. Let $c\mc{M}$ be the category of cosimplicial objects over \mc{M}. 
We refer to \cite{GoJar:simp}, \cite{Hir:loc} or \cite{Hov:model} for the necessary background, in particular for the internal simplicial structure, which is compatible with the Reedy structure, and for latching- and matching objects. Beware of a degree shift between our matching objects and the ones in \cite{GoJar:simp}. The theory of \mc{G}-structures does not require a simplicial structure on \mc{M}, but since we want later on a good theory of $\Tot\co c\mc{M}\to\mc{M}$, we assume from the beginning that the model structure on \mc{M} is simplicial.
We review in subsection \ref{subsection:G-Struktur} the relevant definitions from \cite{Bou:cos} on resolution model structures. In subsection \ref{subsection:natur+spirale} we give a dualized account of the spiral exact sequence from \cite{DKSt:bigraded}. In the last subsection \ref{subsection:Turm} we explain the truncated model structures from \cite{Biedermann:truncated}.

\subsection{The \mc{G}-structure on $c\mc{M}$}
\label{subsection:G-Struktur}
The following definitions are taken from \cite{Bou:cos} who gave the definitive treatment on resolution model structures. 

\begin{Def} \label{injektive Modelle}
Let \mc{M} be a left proper pointed model category. We call a class \mc{G} of objects in \mc{M} a class of {\bf injective models} if the elements of \mc{G} are fibrant and group objects in the homotopy category \ho{\mc{M}} and if \mc{G} is closed under loops. We reserve the letter \mc{G} for such a class.
\end{Def} 

\begin{Def} \label{G-monisch, G-injektiv}
A map $A\toh{i} B$ in \mc{M} is called {\bf \mc{G}-monic} when $\stabhom{B}{G}\toh{i^*}\stabhom{A}{G}$ is surjective for each $G\aus\mc{G}$. 

An object $I$ is called {\bf\boldmath$\mc{G}$\unboldmath-injective} when $\stabhom{B}{I}\toh{i^*}\stabhom{A}{I}$ is surjective for each \mc{G}-monic map $A\toh{i} B$. 

We call a fibration in \mc{M} a {\bf\boldmath$\mc{G}$\unboldmath-injective fibration} if it has the right lifting property with respect to every \mc{G}-monic cofibration.

We say that \ho{\mc{M}} {\bf has enough \boldmath$\mc{G}$\boldmath-injectives} if each object in \ho{\mc{M}} is the source of a \mc{G}-monic map to a \mc{G}-injective target. We say that \mc{G} is {\bf functorial}, if these maps can be chosen functorially.
\end{Def}

\begin{Def} 
A map $X^{\bullet}\to Y^{\bullet}$ in $c\mc{M}$ is called 
\begin{punkt}
   \item
a {\bf\boldmath \mc{G}-equivalence} if the induced maps $\stabhom{Y^{\bullet}}{G}\to\stabhom{X^{\bullet}}{G}$ are weak equivalences of simplicial sets for each $G\aus\mc{G}$.
   \item 
a {\bf\boldmath \mc{G}-cofibration} if it is a Reedy cofibration and the induced maps $\stabhom{Y^{\bullet}}{G}\to\stabhom{X^{\bullet}}{G}$ are fibrations of simplicial sets for each $G\aus\mc{G}$.
   \item
a {\bf\boldmath \mc{G}-fibration} if $f:X^n\to Y^n\times_{M^nY^{\bullet}}M^nX^{\bullet}$ is a \mc{G}-injective fibration for $n\bth 0$.
\end{punkt}
Here $M^nX^\bullet$ is the $n$-th matching object, see \cite[15.2.2.]{Hir:loc}.
These three classes of \mc{G}-equivalences, \mc{G}-cofibrations and \mc{G}-fibrations will be called the {\bf\boldmath \mc{G}-structure on $c\mc{M}$\unboldmath}. We denote it by {\boldmath $c\mc{M}^{\mc{G}}$}.
If \mc{G} is a class of injective models, then by \cite[3.3.]{Bou:cos} the \mc{G}-structure is a simplicial left proper model structure on $c\mc{M}$. 
The simplicial structure is the external one described in appendix \ref{appendix:external}.
\end{Def}

\subsection{Natural homotopy groups and the spiral exact sequence}
\label{subsection:natur+spirale}

It will be important for us to have a different view on the \mc{G}-structure. The reason is, that there are no Postnikov-like truncations with respect to the groups $\pi_s\stabhom{X^\bullet}{G}$. We need a more (co-)homotopical description of the \mc{G}-equivalences. But first we rewrite the above approach and associate groups to a cosimplicial object, which the reader should consider as its cohomology. 
They are contravariant functors on $c\mc{M}$ and depend on two parameters.
In the situation of definition \ref{injektive Modelle} let $X^\bullet$ be an object in $c\mc{M}$ and let $\text{ho}\,\mc{G}$ be the class $\mc{G}$ considered as a full subcategory of $\ho{\mc{M}}=\mc{T}$. Let $\stabhom{\frei}{\frei}$ denote the morphisms in \mc{T}. 
Note that $\stabhom{X^\bullet}{G}$ is a simplicial group.
For every $s\bth 0$ we have a functor
\begin{equation}\begin{split}\begin{array}{ccc}
    \text{ho}\,\mc{G} & \to & \text{groups} \\
       G   & \mapsto & \pi_s\stabhom{X^\bullet}{G} .
    \end{array}\end{split} 
\end{equation}
which takes values in abelian groups for $s\uber 0$. 
Obviously \mc{G}-equivalences are characterized by these groups.
On the other hand we can consider the pointed simplicial set $\Hom{\mc{M}}{X^\bullet}{G}$, where the constant map $X^0\to G$ of the pointed category \mc{M} serves as basepoint. 
Note that if $X^\bullet$ is Reedy cofibrant then this simplicial set is fibrant, and a homotopy group object, since $G$ is so. 
It supplies a functor
\begin{equation}\begin{split}\begin{array}{ccc}
    \mc{G} & \to & \text{fibrant homotopy group objects} \\
       G   & \mapsto & \Hom{\mc{M}}{X^\bullet}{G} . 
             \end{array}\end{split}
\end{equation}
where $\mc{G}$ is considered as a full subcategory of \mc{M}.
Its homotopy should be thought of as the (co-)homotopy of $X^\bullet$. Observe also the equality:
\begin{equation} \label{Hom=map}
    \Hom{\mc{M}}{X^\bullet}{G}= \map^{\rm ext}\hrl X^\bullet, r^0G\hrr ,
\end{equation}
where $r^0G$ denotes the constant cosimplicial object over $G$. Here $\map^{\rm ext}$ denotes the mapping space from the external simplicial structure on $c\mc{M}$ described in appendix \ref{appendix:external}. We will usually drop the superscript.

\begin{Def} \label{natural homotopy groups}
Following \cite{GoHop:moduli} we denote the homotopy groups of these $H$-spaces by
    $$ \naturalpi{s}{X^\bullet}{G} := \pi_s\Hom{\mc{M}}{X^\bullet}{G} = \pi_s\map\hrl X^\bullet,r^0G\hrr $$
for $s\bth 0$ and $G\aus\mc{G}$ and call them the {\bf\boldmath natural homotopy groups of $X^\bullet$ with coefficients in $\mc{G}$}. 
Note that $r^0G$ is Reedy fibrant, so, again, these groups have homotopy meaning if $X^\bullet$ is Reedy cofibrant.
\end{Def}

\begin{Bem} 
Obviously the canonical functor $\mc{M}\to\ho{\mc{M}}$ induces a map $\Hom{\mc{M}}{X^\bullet}{G}\to\stabhom{X^\bullet}{G}$ which in turn induces a natural transformation of functors
\begin{equation}\nonumber 
     \naturalpi{s}{X^\bullet}{G}\to \pi_s\stabhom{X^\bullet}{G} .
\end{equation}
This map is called the {\bf Hurewicz map} and was constructed in \cite[7.1]{DKSt:bigraded}. 
One of the main results is \cite[8.1]{DKSt:bigraded} (also \cite[3.8]{GoHop:moduli}) that this Hurewicz homomorphism for each $G\aus\mc{G}$ fits into a long exact sequence, the so-called {\bf spiral exact sequence}
\begin{align*}
    ...& \to \naturalpi{s-1}{X^\bullet}{\Omega G}\to\naturalpi{s}{X^\bullet}{G}\to\pi_s\stabhom{X^\bullet}{G}\to \naturalpi{s-2}{X^\bullet}{\Omega G} \to ... \\
    ...& \to \pi_{2}\stabhom{X^\bullet}{G}\to \naturalpi{0}{X^\bullet}{\Omega G}\to \naturalpi{1}{X^\bullet}{G}\to \pi_1\stabhom{X^\bullet}{G} \to 0,
\end{align*}
where $\Omega$ is the loop space functor on \mc{M}, plus an isomorphism 
    $$ \naturalpi{0}{X^\bullet}{G}\cong\pi_0\stabhom{X^\bullet}{G}.  $$
In the construction of the exact sequence we rely on the external simplicial structure, but not on a simplicial structure on \mc{M}. 

As explained in \cite[8.3.]{DKSt:bigraded} or \cite[(3.1)]{GoHop:moduli} these long exact sequences can be spliced together to give an exact couple and an associated spectral sequence
\begin{equation} \label{Spiralspektral}
      \pi_p\stabhom{X^\bullet}{\Omega^q G}\ \Longrightarrow\ \colim_k\naturalpi{k}{X^\bullet}{\Omega^{p+q-k}G} .
\end{equation}
\end{Bem}

\begin{lemma} \label{naturale G-Äquivalenzen}\label{naturale G-Aequivalenzen} 
A map $X^\bullet\to Y^\bullet$ is an \mc{G}-equivalence if and only if it induces iso\-mor\-phisms
    $$ \emph{\naturalpi{s}{\wt{Y}^\bullet}{G}}\to\emph{\naturalpi{s}{\wt{X}^\bullet}{G}} $$
for all $s\bth 0$ and all $G\in\mc{G}$ and some Reedy cofibrant approximation $\wt{X}^\bullet\to\wt{Y}^\bullet$.
\end{lemma}

\begin{beweis}
This follows immediately from the spiral exact sequence by simultaneous induction over the whole class \mc{G} and the five lemma. Remember that \mc{G} is closed under loops by assumption.
\end{beweis}

\subsection{Truncated \mc{G}-structures}
\label{subsection:Turm}
Now we will study the homotopy categories associated to the truncated structures. Observe, that there is a natural isomorphism
    $$ \map\hrl\sk_{n+1}X^\bullet, r^0G\hrr\cong\cosk_{n+1}\map\hrl X^\bullet, r^0G\hrr .$$
If $X^\bullet$ is Reedy cofibrant then $\map\hrl X^\bullet, r^0G\hrr$ is fibrant. Note also, that for a Kan complex $W$ the space $\cosk_{n+1}W$ is a model for the $n$-th Postnikov section.

\begin{Def} \label{Def. n-G-Struktur}
In \cite[theorem 3.5.]{Biedermann:truncated} we prove the existence of a left proper simplicial model structure on $c\mc{M}$, called the {\bf\boldmath $n$-\mc{G}-structure}, whose equivalences are given by maps $X^\bullet\to Y^\bullet$ such that for every $G\in\mc{G}$ and all $0\sth s\sth n$ the induced maps
    $$ \naturalpi{s}{\wt{Y}^\bullet}{G}\to\naturalpi{s}{\wt{X}^\bullet}{G} $$
are isomorphisms, where $\wt{X}^\bullet\to\wt{Y}^\bullet$ is a cofibrant approximation to $X^\bullet\to Y^\bullet$. These maps are called {\bf\boldmath $n$-\mc{G}-equivalences}. An {\bf\boldmath $n$-$\mc{G}$-cofibration} is a map $X^\bullet\to Y^\bullet$ which is a \mc{G}-cofibration such that for every $G\in\mc{G}$ and all $s\uber n$ the induced maps    
    $$ \naturalpi{s}{\wt{Y}^\bullet}{G}\to\naturalpi{s}{\wt{X}^\bullet}{G} $$
are isomorphisms.
The fibration are the \mc{G}-fibrations.
\end{Def}

Now we are going to determine the $n$-\mc{G}-cofibrant objects. 
\begin{Bem} \label{n-G-kofibrante Objekte}
Remember that cofibrant objects in the \mc{G}-structure coincide with the Reedy cofibrant ones.
\begin{punkt}
      \item An object $A^\bullet$ in $c\mc{M}$ is $n$-\mc{G}-cofibrant if and only if it is Reedy cofibrant and $\naturalpi{s}{A^\bullet}{G}=0$ for all $G\aus\mc{G}$ and $s\uber n$. 
      \item An $n$-\mc{G}-cofibrant approximation functor is given by $Q=\sk_{n+1}\wt{\frei}$.
      \item
On $n$-\mc{G}-cofibrant objects the $n$-\mc{G}-structure and the \mc{G}-structure coincide. 
\end{punkt}
\end{Bem}

\begin{Def} \label{Ko-Postnikov-Turm}
Let $X^\bullet$ be an object in $c\mc{M}$. The skeletal filtration of a Reedy cofibrant approximation to $X^\bullet$ consists of $n$-\mc{G}-cofibrant approximations $X^\bullet_n$ to $X^\bullet$ for the various $n$, and these assemble into a sequence
    $$ X^\bullet_0\to X^\bullet_1\to X^\bullet_2\to ... \to X^\bullet $$
which captures higher and higher natural homotopy groups. So this can be viewed as a {\bf Postnikov cotower} for $X^\bullet$. 
\end{Def}

\begin{Def} \label{Turm von abgeschnittenen Homotopiekategorien}
As explained in \cite[3.13.]{Biedermann:truncated} the functor $\id\co c\mc{M}^{\mc{G}}\to c\mc{M}^{n-\mc{G}}$ is a right Quillen functor, whose left adjoint is given by $Q_n=\sk_{n+1}\wt{\frei}$. We have an induced pair of adjoint derived functors:
    $$ LQ_n\cong L(\id)\co\ho{c\mc{M}^{n-\mc{G}}}\leftrightarrows\ho{c\mc{M}^{\mc{G}}}:\!R(\id), $$
where $LQ_n\cong L(\id)$ is an embedding of a full subcategory. In the same way, we can view $\id\co c\mc{M}^{(n+1)-\mc{G}}\to c\mc{M}^{n-\mc{G}}$ as a right Quillen functor and $L(\id)\co\ho{c\mc{M}^{n-\mc{G}}}\to\ho{c\mc{M}^{(n+1)-\mc{G}}}$ is again an embedding of a full subcategory.
The tower of categories
    $$ ... \to \ho{c\mc{M}^{(n+1)-\mc{G}}} \toh{\sigma_n} \ho{c\mc{M}^{n-\mc{G}}}\to ... \to \ho{c\mc{M}^{1-\mc{G}}}\toh{\sigma_0} \ho{c\mc{M}^{0-\mc{G}}} $$
can be identified as a tower of full subcategories of $\ho{c\mc{M}^{\mc{G}}}$ given by coreflections.

We can characterize the objects in $\ho{c\mc{M}^{n-\mc{G}}}$ viewed as a subcategory of $\ho{c\mc{M}^{\mc{G}}}$ by their natural homotopy groups. An object $X^\bullet$ is in the image of $\ho{c\mc{M}^{n-\mc{G}}}$ if and only if it is \mc{G}-equivalent to its $n$-\mc{G}-cofibrant replacement, i.e. if we have:  
\begin{equation*}
    \naturalpi{s}{\wt{X}^\bullet}{G}=0 \hbox{ for } s\uber n
\end{equation*}
for a Reedy cofibrant replacement $\wt{X}^\bullet\to X^\bullet$.
We have to relate all this to $\ho{\mc{M}}$ by the following statement, whose analogue for the Reedy structure is well known. The lemma is cited from \cite[Prop. 8.1.]{Bou:cos}.
\end{Def}

\begin{lemma} \label{Tot-Delta-Quillen-Paar}
The functors 
    $$\xymatrix{ \mc{M} \ar@<2pt>[rr]^-{\frei\otimes^{\rm pro}\Delta^\bullet} & & c\mc{M}^{\mc{G}} \ar@<2pt>[ll]^-{\Tot} } .$$
form a Quillen pair.
\end{lemma}

\begin{Bem} \label{Dasselbe Quillenpaar}
The natural transformation $\frei\otimes^{\rm pro}\Delta^\bullet\to r^0$ gives a Reedy cofibrant replacement by \cite[16.1.4.]{Hir:loc} and hence a \mc{G}-cofibrant replacement. It follows that both induce the same left derived functor:
    $$\xymatrix{ \ho{\mc{M}} \ar@<2pt>[rrr]^-{Lr^0=\,\frei\otimes^{\rm pro}\Delta^\bullet} &&& \ho{c\mc{M}^{\mc{G}}} \ar@<2pt>[lll]^-{R\Tot} } .$$
We can look at the composition
    $$ \mc{M}\toh{r^0} c\mc{M}^{\mc{G}} \toh{\id} c\mc{M}^{n-\mc{G}}  $$
and the composition of induced derived functors:
\begin{equation}\label{theta_n} 
      \xymatrix{\ho{\mc{M}}\ar[r]^-{Lr^0} & \ho{c\mc{M}^{\mc{G}}} \ar[r]^-{R(\id)} & \ho{c\mc{M}^{n-\mc{G}}} }.  
\end{equation}
\end{Bem}

\begin{Def} 
We will denote the composition \Ref{theta_n} of functors by \boldmath$\theta_n$\unboldmath. We arrive at the following diagram:
     $$\xymatrix@=20pt{ && \mc{T}=\ho{\mc{M}} \ar[d]^{\theta_n}\ar[dl]_{\theta_{n+1}}\ar[drr]^{\theta_0} && \\
        ... \ar[r] & \ho{c\mc{M}^{(n+1)-\mc{G}}} \ar[r]_-{\sigma_n} & \ho{c\mc{M}^{n-\mc{G}}}\ar[r] &  ... \ar[r]_-{\sigma_0} & \ho{c\mc{M}^{0-\mc{G}}}  } $$
This diagram is a $2$-commuting diagram of functors. For details on $2$-commutativity we refer to \cite{Hov:model}.
$2$-commutativity is provided by the relation $QQ\simeq Q$.
We call this {\bf the tower of truncated homotopy categories} associated to \mc{M} and \mc{G}. 
\end{Def}

\section{Homological functors} 
\label{section:homologische Funktoren}

This section begins in \ref{subsection:homologische Funktoren mit genügend Injektiven} with the introduction of the technical conditions, that the homological functors we consider have to satisfy. In particular, the general assumptions on $F$ are summarized in \ref{Annahmen fuer F}. In subsection \ref{subsection:F-injektive Struktur} we derive the resolution model structures, which are relevant for the obstruction calculus, by applying the general machinery from section \ref{section:E2-Strukturen}. We study the associated homotopy category in \ref{subsection:F-Homotopiekategorie}.

\subsection{Homological functors with enough injectives}
\label{subsection:homologische Funktoren mit genügend Injektiven}

\begin{Def} \label{Einhängungen}
From now on let \boldmath$\mc{T}$\unboldmath\ always denote a triangulated category. The set of morphisms for $X$ and $Y$ in \mc{T} will be denoted by \boldmath$\stabhom{X}{Y}$\unboldmath. 
The shift functor or {\bf suspension} of \mc{T} will be denoted by \boldmath$\Sigma$\unboldmath. It is, of course, an equivalence of categories.

Let \boldmath$\mc{A}$\unboldmath\ always be a {\bf graded abelian category}, which means, we require that \mc{A} possesses a {\bf shift functor} denoted by $\hel 1\her$ which is an equivalence of categories. Let $\hel n\her$ denote the $n$-fold iteration of $\hel 1\her$.
\end{Def}

\begin{Def} \label{homologische Funktoren}
Let $F_*:\mc{T}\to\mc{A}$ be a covariant functor, where the star stands for the grading of \mc{A}. We say that $F_*$ is {\bf homological}, if it satisfies the following conditions: 
\begin{punkt}
    \item $F_*$ is a graded functor, in other words, it commutes with suspensions, so there are a natural equivalences 
          $$ F_*\Sigma X\cong \hrl F_*X\hrr\hel 1\her =: F_{*-1}X,$$
which are part of the structure.
    \item $F_*$ is additive saying that it commutes with arbitrary coproducts.
    \item $F_*$ converts distinguished triangles into long exact sequences.
\end{punkt}
\end{Def}

\begin{Bem} 
Later we will assume that \mc{T} has an underlying model category \mc{M}.
The suspension functor $\Sigma$ here is internal to the model structure on \mc{M}. The reader should be aware that this has nothing to do with the external construction $\Sigma_{\rm ext}$ from \ref{aeussere Einhaengung} which is derived from the external simplicial structure on the cosimplicial objects $c\mc{M}$ over \mc{M}. 
\end{Bem}

\begin{Def} \label{F-Lokalität}
We say that {\bf\boldmath $F:\mc{T}\to\mc{A}$ detects isomorphisms\unboldmath} or equivalently that {\bf\boldmath \mc{T} is $F$-local\unboldmath} if a map $X\to Y$ in \mc{T} is an isomorphism if and only if the induced map $F_*X\to F_*Y$ in \mc{A} is an isomorphism.
\end{Def}

\begin{Def}  
\label{Eilenberg-MacLane} 
Let $F_*$:$\mc{T} \to \mc{A}$ be a homological functor from a triangulated category to a graded abelian category and let $I$ be an injective object in \mc{A}.
Consider the following functor:
\begin{equation}
      X \mapsto \Hom{\mc{A}}{F_*X}{I} \nonumber
\end{equation} 
We require that this functor is representable by an object $E\hrl I\hrr$ of \mc{T}. If the canonical morphism $F_*E\hrl I\hrr\to I$ induced by
\begin{equation}
      {\rm id}_{E(I)}\in\stabhom{E\hrl I\hrr}{E\hrl I\hrr} \cong \Hom{\mc{A}}{F_*E\hrl I\hrr}{I} \nonumber
\end{equation}
is an isomorphism, then we call $E\hrl I\hrr$ an {\bf\boldmath $\hrl F,I\hrr$-Eilenberg-MacLane object\unboldmath}. 
Usually we will just say that $E\hrl I\hrr$ is {\bf\boldmath $F$-injective\unboldmath}. 
\end{Def}

\begin{Def} \label{genügend Injektive}
We will say that the functor {\bf $F_*$ possesses enough injectives}, if every object in \mc{T} admits a morphism to an $F$-injective object that induces a monomorphism in \mc{A}.
\end{Def}

\begin{Bem} \label{Anwendungen}
There are a lot of examples of functors with enough injectives that detect isomorphisms. Every topologically flat ring spectrum $E$, where $E_*E$ is commutative, induces a homological functor
     $$ E_*:\ho{\text{Spectra}}\to E_*E{\rm -comod} $$
from the stable homotopy category of spectra to the category of $E_*E$-comodules. The phrase ring spectrum is to be interpreted here in the most naive sense, a monoid object in the homotopy category. The notion of topological flatness was earlier called the Adams(-Atiyah) condition. It ensures, see \cite[Thm. 1.5.]{Dev:brown-comenetz}, that $E_*E(I)\cong I$ as $E_*E$-comodules. 
Hence $E_*$ possesses enough injectives, since $F$-injective objects exist by Brown representability. See also \cite{Hov:algebroid}. If we perform Bousfield localization and consider it to be a functor
     $$ E_*:\ho{\text{Spectra}}_E\to E_*E{\rm -comod} $$
from the $E$-local category, it also detects isomorphism. From this data we can construct a spectral sequence, which is known as the $E$-based modified Adams-spectral sequence, see \ref{ASS}. 
\end{Bem}

\begin{Bem} \label{F-Injektive und Retrakte}
There are two convenient facts about finding $F$-injectives that are derived from \cite[2.1.1.]{Fra:uni} or \cite[Thm. 1.5.]{Dev:brown-comenetz}:
\begin{punkt}
      \item
If $F$ detects isomorphisms, then every representing object $X$ in \mc{T}, whose image $F_*X\cong I$ is injective in \mc{A}, is an $\hrl F,I\hrr$-Eilenberg-MacLane object.
      \item
Retracts of $F$-injective objects are again $F$-injective.
\end{punkt}
\end{Bem}

\begin{Bem} \label{E(I) ist funktoriell}
If $F_*$ is a homological functor that possesses enough injectives and detects isomorphism, then there is the following observation taken from \cite[2.1. Lemma 1]{Fra:uni}: Given another representing object $\wt{E}\hrl I\hrr$, there is a unique morphism
     $$ \wt{E}\hrl I\hrr\to E\hrl I\hrr $$
in \mc{T} lifting the identity of $I$, and this is an isomorphism.
This can be reformulated in the following way:
Let $\mc{T}_{F{\rm -inj}}$ denote the full subcategory of \mc{T} consisting of the $F$-injective objects. Let $\mc{A}_{\rm inj}$ denote the full subcategory of \mc{A} consisting of the injective objects. Then the functor $F$ induces an equivalence $\mc{T}_{F{\rm -inj}}\to\mc{A}_{\rm inj}$.
\end{Bem}

The following assumptions will be valid for the rest of the article.
\begin{Ann} \label{Annahmen für F}\label{Annahmen fuer F}
From now on let \mc{T} be the homotopy category of a simplicial left proper stable model category \mc{M} and let \mc{A} be an abelian category with enough injectives. 
Let $F_*:\mc{T}\to\mc{A}$ be a homological functor with enough $F$-injectives as explained in \ref{homologische Funktoren} and \ref{genügend Injektive}, which detects isomorphisms \ref{F-Lokalität}. 
We will call the composition $\mc{M}\to\mc{T}\to\mc{A}$ of $F_*$ with the canonical functor from \mc{M} to its homotopy category also $F_*$. By applying it levelwise we can prolong it to a functor $c\mc{M}\to c\mc{A}$ that we will again call $F_*$.
\end{Ann}

\subsection{The $F$-injective structure and its truncations}
\label{subsection:F-injektive Struktur}

\begin{Def} \label{G=F-Inj}
We take as our class of injective models \mc{G} the class of all $F$-injective objects in \mc{M} which were defined in \ref{Eilenberg-MacLane} . This class \mc{G} will be fixed for the rest of this work. 
We denote our special choice by 
\begin{center}
    $\{F\text{-injectives}\}=:$ {\boldmath$\{F$\unboldmath}{\bf -Inj}{\boldmath$\}$\unboldmath}. 
\end{center} 
We will call the involved classes of maps {\bf\boldmath$F$\unboldmath -injective equivalences}, {\bf\boldmath$F$\unboldmath -injec\-tive fibrations} and {\bf cofibrations}. Sometimes we will abbreviate even this and simply say $F$-equivalent or $F$-fibrant and so on.
We will call this model structure on $c\mc{M}$ the {\bf \boldmath{$F$}-injective model structure}.
The truncated model structure from \ref{Def. n-G-Struktur} will be called {\bf\boldmath $n$-$F$-injective structure\unboldmath} or just $n$-$F$-structure. 
We will denote them by $c\mc{M}^F$ and $c\mc{M}^{n-F}$.
\end{Def}

To conclude that the choice of $\mc{G}=\{F\text{-injectives}\}$ is really admissible, we observe first of all, that \mc{M} is stable, so all objects are homotopy group objects. Next it is easy to check from the definitions, that the $F$-injectives are closed under equivalences and (de-)suspensions. Finally it will follow, that there are enough \mc{G}-injectives, from the assumption, that $F$ has enough injectives \ref{genügend Injektive} and the following two consistency checks, which are easy to prove.
A map will be called {\bf \boldmath{$F$}-monic} if it is $\{F\text{-Inj}\}$-monic. 

\begin{lemma} \label{F-monisch}
A map is $F$-monic if and only if its image is a monomorphism. 
\end{lemma}

\begin{lemma} \label{Saturiertheit der F-Injektiven}
The two classes $\{F{\rm -Inj}\}$ and $\{\{F{\rm -Inj}\}\text{-{\rm injectives}}\}$ coincide.
\end{lemma}

As a consequence of \cite{Biedermann:truncated} -- see also \ref{Def. n-G-Struktur} -- we have the following model structures at hand, where right properness is proved in \ref{G-Struktur ist eigentlich}.
\begin{satz} \label{(abgeschnittene) F-Strukturen}
Let \mc{M} be a pointed simplicial left proper stable model category und set $\ho{\mc{M}}=:\mc{T}$ and let \mc{A} be an abelian category. Let $F:\mc{T}\to\mc{A}$ be a homological functor that possesses enough injectives and that detects isomorphisms. On $c\mc{M}$ there is a pointed simplicial proper model structure given by the $\emph{\hrl} n$-$\emph{\hrr}F$-injective equivalences, the $\emph{\hrl} n$-$\emph{\hrr}F$-injective cofibrations and the $\emph{\hrl} n$-$\emph{\hrr}F$-injective fibrations.
The simplicial structure is always the external one.
\end{satz}

In fact $\ho{c\mc{M}}^F$ behaves like the category of non-negative cochain complexes inside the full derived category of an abelian category with enough injectives.
See the discussion in subsection \ref{subsection:F-Homotopiekategorie}.

\begin{Bem} \label{injektive Struktur auf cA} 
If we view \mc{A} as a discrete model category, we can equip the category $c\mc{A}$ of cosimplicial objects over \mc{A} with the {\bf\boldmath \mc{I}-structure}, where the class \mc{I} of injective objects in \mc{A} is taken as a class of injective models. It follows from \cite[4.4]{Bou:cos} that this model structure corresponds to the classical model structure from \cite{Qui:htp} for the nonnegative cochain complexes ${\rm CoCh}^{\bth 0}\hrl\mc{A}\hrr$ via the Dold-Kan correspondence. So in ${\rm CoCh}^{\bth 0}\hrl\mc{A}\hrr$ we have:
The \mc{I}-equivalences are the cohomology equivalences, the \mc{I}-co\-fi\-brations are the maps that are monomorphisms in positive degrees, and the \mc{I}-fi\-bra\-tions are those that are (split) surjective with injective kernel in all degrees. 
The fibrant objects are the degreewise injective objects, while all objects are cofibrant.
\end{Bem}

We will now list characterizations of $F$-injective equivalences, $F$-injective cofibrations and $F$-injective fibrations and their truncated analogues. 
In the next statements let $\pi^sA^\bullet\cong H^sNA^\bullet$ be the usual thing with many names, e.g. the cohomology of the normalized cochain complex $NA^\bullet$.

\begin{kor} \label{F-injektive Äquivalenzen}
A map $X^{\bullet}\to Y^{\bullet}$ in $c\mc{M}$ is an $F$-injective equivalence if and only if the induced maps
    $$ H^sNF_*X^{\bullet}\to H^sNF_*Y^{\bullet}$$ 
are isomorphisms for all $s\bth 0$, i.e. it induces a quasi-isomorphism $NF_*X^\bullet\to NF_*Y^\bullet$.
\end{kor}

\begin{beweis}
We have for $F$-injective $G$ the isomorphisms
    $$ \pi_s\stabhom{X^\bullet}{G}\cong H^sN\Hom{\mc{A}}{F_*X^\bullet}{F_*G}\cong \Hom{\mc{A}}{H^sNF_*X^\bullet}{F_*G}$$
Then the lemma follows from the fact mentioned in \ref{E(I) ist funktoriell}, that if $G$ runs through all $F$-injectives then $F_*G$ ranges over all injectives in \mc{A}.
\end{beweis}

\begin{Bem} \label{nichts Offensichtliches}
There is no obvious way to characterize $n$-$F$-equivalences in terms of $\pi^sF_*\hrl\frei\hrr$ like the $F$-equivalences. The induction used to prove \ref{naturale G-Aequivalenzen} crawling up the spiral exact sequence does not yield anything useful if it stops at some finite stage. So we do not offer another description of them as the one given in \ref{Def. n-G-Struktur}.
Of course we have, that on $n$-$F$-cofibrant objects $n$-$F$-equivalence and $F$-equivalence agree.
\end{Bem}

\begin{lemma} \label{F-Kofaserungen}
A map $i:X^{\bullet}\to Y^{\bullet}$ is an $F$-injective cofibration if and only if it is a Reedy-cofibration that induces monomorphisms
    $$ N^k FX^{\bullet} \to N^k FY^{\bullet} $$
for all $k\bth 1$.
\end{lemma}

\begin{beweis}
The map $i$ is an $F$-cofibration if and only if it is a Reedy cofibration and the induced map
     $$ \stabhom{Y^\bullet}{G}\to\stabhom{X^\bullet}{G} $$
is a fibration of simplicial sets for all $G\aus\mc{G}$. The result now follows from the fact that a map of simplicial abelian groups is a fibration if and only if it induces a surjection of the normalizations in positive degrees.
\end{beweis}

\begin{lemma} \label{lange exakte Homologiesequenz}
Let $X^\bullet\to Y^\bullet$ be an $F$-cofibration with cofiber $C^\bullet$ that induces a mo\-no\-mor\-phism $N^0F_*X^\bullet\to N^0F_*Y^\bullet$. Then there is a long exact sequence 
\begin{align*}
     0\to H^0NF_*X^\bullet\to H^0NF_*Y^\bullet\to H^0NF_*C^\bullet\to H^1NF_*X^\bullet\to ... \\
     ...\to H^sNF_*X^\bullet\to H^sNF_*Y^\bullet\to H^sNF_*C^\bullet\to H^{s+1}NF_*X^\bullet\to ...  
\end{align*}
\end{lemma}

\begin{beweis}
This can be proved by 
\ref{F-Kofaserungen}.
\end{beweis}

We need a little bit more care to describe $F$-injective fibrations. First of all we remind the reader that $F$-injective fibrations and $n$-$F$-injective fibrations coincide. By definition a map $X^\bullet\to Y^\bullet$ is an $F$-injective fibration if and only if all the maps $X^s\to M^sX^\bullet\times_{M^sY^\bullet}Y^s$ for $s\bth 0$ are \mc{G}-injective fibrations in \mc{M} in the sense of definition \ref{G-monisch, G-injektiv} with $\mc{G}=\{F\text{-Inj}\}$. Thus we describe $\{F\text{-Inj}\}$-injective fibrations in \mc{M}.
\begin{lemma} \label{F-injektive Faserungen in M}
A map in \mc{M} is an \{$F${\rm -Inj}\}-injective fibration if and only if it is a fibration with $F$-injective fiber and that induces an epimorphism under $F$.
\end{lemma}

\begin{beweis}
By \cite[3.10.]{Bou:cos} a map $X\to Y$ in \mc{M} is a \mc{G}-injective fibration if and only if it is a retract of a \mc{G}-cofree map $X'\to Y'$. A \mc{G}-cofree map is a map that can be expressed as a composition $X'\to Y'\times E\to Y'$, where $X'\to Y'\times E$ is a trivial fibration in \mc{M}, $Y'\times E\to Y'$ is the projection onto $Y'$ and $E$ is \mc{G}-injective. 

The assertion is true for $\{F\text{-Inj}\}$-co\-free maps. Here we use the fact that \mc{T} is $F$-local (see \ref{F-Lokalität}), so weak equivalences in \mc{M} induce isomorphisms under $F$.
 But the claim is also true for retracts. This is obvious for surjectivity. The fiber condition follows from \ref{F-Injektive und Retrakte}, since the fiber of $X\to Y$ is a retract of the fiber of $X'\to Y'$ which is weakly equivalent to $E$ and therefore itself $F$-injective.

Conversely let $X\to Y$ be a fibration that has an $F$-injective fiber $E$ and that induces a surjection under $F$. \mc{M} is stable, hence we get a long exact $F$-sequence for $X\to Y$ and it follows $X\simeq E\times Y$. We deduce that $X\to Y$ has the right lifting property with respect to every $\{F\text{-Inj}\}$-monic cofibration. So it is an $\{F\text{-Inj}\}$-injective fibration. 
\end{beweis}

\begin{kor} \label{F-Faserungen}
A Reedy fibration $X^\bullet\to Y^\bullet$ between $F$-fibrant objects is an $F$-fibration if and only if for $s\bth 0$ the induced maps $N^sF_*X^\bullet\to N^sF_*Y^\bullet$ are surjective with injective kernel. In other words, $X^\bullet\to Y^\bullet$ is an $F$-fibration if and only if $F_*X^\bullet\to F_*Y^\bullet$ is an \mc{I}-fibration.
\end{kor}

\begin{beweis}
This follows from \ref{F-injektive Faserungen in M}, since for an $F$-fibrant $Y^\bullet$ and all $s\bth 0$ we have an isomorphism
    $$ F_*\hrl Y^s\times_{M^sY^\bullet}M^sX^\bullet\hrr\cong F_*Y^s\times_{F_*M^sY^\bullet}F_*M^sX^\bullet.$$
\end{beweis}

\begin{lemma} \label{F-Homotopiepullbacks}
The functor $F_*:c\mc{M}\to c\mc{A}$ maps $F$-homotopy pullbacks to \mc{I}-homotopy pullbacks.
\end{lemma}

\begin{beweis}
It is sufficient to prove that $F_*$ preserves pullbacks after $F$-fibrant replacement. Let $X^\bullet\to Z^\bullet\leftarrow Y^\bullet$ be $F$-fibrations between $F$-fibrant objects. It follows by \cite[5.3.]{Bou:cos} that all maps $X^s\to Z^s\leftarrow Y^s$ are $F$-injective fibrations in \mc{M}, in particular they are fibrations and induce surjections under $F_*$. The pullback square
\diagr{ X^s\times_{Z^s}Y^s \ar[r]\ar[d] & X^s \ar[d] \\ Y^s \ar[r]& Z^s}
is also a homotopy pullback square in \mc{M}, and hence a homotopy pushout. 
The long exact $F$-sequence resulting from this collapses to the short exact sequences
    $$ 0\to F_*\hrl X^s\times_{Z^s}Y^s\hrr \to F_*X^s\oplus F_*Y^s \to F_*Z^s \to 0 ,$$
which proves the lemma.
\end{beweis}

\begin{kor} \label{G-Struktur ist eigentlich}
The $(n$-$)F$-structure for $0\sth n\sth\infty$ is proper.
\end{kor}

\begin{beweis}
We only need to prove right properness. For $n=\infty$ this follows from the characterization of $F$-equivalences and $F$-fibrations in \ref{F-injektive Äquivalenzen} and \ref{F-Faserungen} and from \ref{F-Homotopiepullbacks}. This passes down to smaller $n$ by theorem \cite[3.5.]{Biedermann:truncated}.
\end{beweis}

\subsection{The $F$-injective homotopy category}
\label{subsection:F-Homotopiekategorie}

The category $c\mc{M}$ equipped with the $F$-structure behaves very much like the full subcategory ${\rm CoCh}^{\bth 0}\hrl\mc{A}\hrr$ of nonnegative cochain complexes within the derived category $D\hrl\mc{A}\hrr$. This is displayed by the statements \ref{fast stabil} and \ref{noch mehr stabil}. We are going to need a dual version of the functor $\ol{W}:s\text{Ab}\to s\text{Ab}$ which is sometimes called the Eilenberg-MacLane functor or the Kan suspension.

\begin{Def} \label{Definition von W und W quer}
Let \mc{N} be a pointed model category.
We define a functor $W:c\mc{N}\to c\mc{N}$. Let $X^\bullet$ be a cosimplicial object. Let \boldmath$W X^\bullet$\unboldmath\ be defined by the following equations:
   $$ \hrl W X^\bullet\hrr^s := \prod_{i=0}^s X^i $$
The structural maps of a cosimplicial object are constructed by the process dual to the one described in \cite[p. 192]{GoJar:simp}. 
There is a map $W X^\bullet\to X^\bullet$ given by projection
     $$ \prod_{i=0}^sX^i \to X^s. $$
Let \boldmath$\ol{W}X^\bullet$\unboldmath\ be the fiber of $W X^\bullet\to X^\bullet$.
\end{Def}

\begin{Bem} \label{Eigenschaften von W und W-quer}
Let $X^\bullet$ be in $c\mc{M}$.  
The map $W X^\bullet\to X^\bullet$ is a Reedy-fibration if and only if every $X^\bullet$ is Reedy fibrant. 
It is a \mc{G}-fibration for some general \mc{G} if and only if in addition all $X^s$ are \mc{G}-injective. In both cases $\ol{W}X^\bullet$ has homotopy meaning, see \ref{mehr Eigenschaften von W und W-quer}.
\end{Bem}

\begin{lemma} 
If we take $\mc{G}=\{F{\rm -Inj}\}$ then $WX^\bullet$ is $F$-equivalent to $*$.
\end{lemma}

\begin{beweis}
Since $FWX^\bullet\cong WFX^\bullet$ it suffices by \ref{F-injektive Äquivalenzen} to show that $WA^\bullet$ is \mc{I}-equivalent to $*$ for arbitrary $A^\bullet$ in $c\mc{A}$. This follows by dualizing \cite[III.5.]{GoJar:simp}. 
\end{beweis}

\begin{Bem} \label{mehr Eigenschaften von W und W-quer}
Hence $\ol{W}X^\bullet$ is another different model for the loop object $\Omega_{\rm ext}X^\bullet$.
If $A^\bullet$ is in $c\mc{A}$ this object can also be obtained in the following way:
\diagr{ c\mc{A} \ar[r]^{\ol{W}}\ar[d]_N & c\mc{A}  \\
        {\rm CoCh}^{\bth 0}\hrl\mc{A}\hrr \ar[r]_{[1]_{\rm ext}} & {\rm CoCh}^{\bth 0}\hrl\mc{A}\hrr \ar[u]_\Gamma}
where $\hrl A^*\hel 1\her_{\rm ext}\hrr^s=A^{s+1}$ is the external shift functor of cochain complexes (which should not be confused with the internal shift $\hel 1\her$ from \ref{Einhängungen}), $N$ is normalization and $\Gamma$ is the Dold-Kan functor. In particular if $A^\bullet$ is in $c\mc{A}$ we have:
\begin{equation} \label{Verschiebung}
       H^s N\ol{W}A^\bullet = \left\{ \begin{array}{cl}
                                  0          &, \hbox{ for } s=0 \\
                           H^{s-1}NA^\bullet &, \hbox{ for } s\bth 1 
                                          \end{array} 
                                  \right.   
\end{equation}
For every $F$-fibrant $X^\bullet$ we get a map
\begin{equation} \label{natürliche Transformation}    
    \Sigma_{\rm ext}\ol{W}X^\bullet \to X^\bullet 
\end{equation}
in $c\mc{M}$ which descends to a natural transformation $\Sigma_{\rm ext}\Omega_{\rm ext} \to \text{Id}$ of endofunctors of $\ho{c\mc{M}^{\mc{G}}}$.
\end{Bem}

\begin{lemma} \label{Omega und W-quer}
For every $F$-fibrant object $X^\bullet$ in $c\mc{M}$ the map $\Sigma_{\rm ext}\ol{W}X^\bullet \to X^\bullet$ is an $F$-equivalence. 
\end{lemma}

\begin{beweis}
We note that $F\Sigma_{\rm ext}\ol{W}X^\bullet=\Sigma_{\rm ext}\ol{W}FX^\bullet$ because $F$ is applied levelwise and commutes with finite products. Now the fact follows from \ref{F-injektive Äquivalenzen} and \Ref{Verschiebung}.
\end{beweis}

\begin{kor} \label{fast stabil}
The map \emph{\Ref{natürliche Transformation}} induces a natural equivalence $\Sigma_{\rm ext}\Omega_{\rm ext}\cong {\rm Id}$ of en\-do\-func\-tors of $\ho{c\mc{M}^{F}}$.
\end{kor}

\begin{beweis}
This follows from \ref{Eigenschaften von W und W-quer}, \ref{mehr Eigenschaften von W und W-quer} and \ref{Omega und W-quer}.
\end{beweis}

Furthermore we have an isomorphism $\Omega_{\rm ext}\Sigma_{\rm ext}X^\bullet\cong X^\bullet$ in $\ho{c\mc{M}^{\mc{G}}}$ as long as the objects in question are ``connected''.

\begin{lemma} \label{noch mehr stabil}
Let $X^\bullet$ be an $F$-fibrant object such that $\pi^0F_*X^\bullet=0$. Then the canonical map $X^\bullet\to\ol{W}\Sigma_{\rm ext}X^\bullet$ is an $F$-equivalence. 
\end{lemma}

\begin{beweis}
The condition $\pi^0F_*X^\bullet=0$ is equivalent to $0=\pi_0\stabhom{X^\bullet}{G}\cong\naturalpi{0}{X^\bullet}{G}$.
Hence the map $X^\bullet\to\ol{W}\Sigma_{\rm ext}X^\bullet$ induces isomorphisms on $H^sNF_*\hrl\frei\hrr$ for all $s\bth 0$, so it is an $F$-equivalence. 
\end{beweis}

\begin{Bem} \label{Biprodukte}
In a stable model category \mc{M} finite products and finite coproducts are weakly equivalent. It follows that for the Reedy structure and in particular for every \mc{G}-structure on $c\mc{M}$ and their truncated versions finite products and coproducts are weakly equivalent. 
\end{Bem}

\begin{kor} \label{Additivität}
For every $0\sth n\sth\infty$ the category $\ho{c\mc{M}^{n-\mc{G}}}$ is additive and the functors $\sigma_n:\ho{c\mc{M}^{(n+1)-F}}\to\ho{c\mc{M}^{n-F}}$ and $\theta_n:\mc{T}\to\ho{c\mc{M}^{n-F}}$ are additive.
\end{kor}

\begin{beweis}
By \ref{fast stabil} every object in $\ho{c\mc{M}^{n-F}}$ for $0\sth n\sth\infty$ is isomorphic to a double suspension, hence every object is an abelian cogroup object in the homotopy category. Both functors $\sigma_n$ and $\theta_n$ commute with $\Sigma_{\rm ext}$.
\end{beweis}

\begin{Def} 
We will denote the biproduct of a pair of objects $X^\bullet$ and $Y^\bullet$ in $\ho{c\mc{M}^{n-F}}$ for $0\sth n\sth\infty$ by \boldmath$X^\bullet\oplus Y^\bullet$\unboldmath.
\end{Def}

\section{The realization problem}

In the first subsection \ref{subsection:Hindernisse} we do the hard work and construct the obstruction calculus. The main theorems are \ref{Ext^(n+2,n) und Objekte}, \ref{Ext^(n+1,n) und Objekte}, \ref{Ext^(n+1,n) und Morphismen}, \ref{Derivation} and \ref{Ext^(n,n) und Morphismen}. 
In \ref{subsection:Turm von Interpolationskategorien} we define our interpolation categories for a homological functor $F$ an in \ref{subsection:SS} we describe the spectral sequences that play a role in the realization problem.

\subsection{Realizations and obstruction calculus}
\label{subsection:Hindernisse}

In this subsection we develop an obstruction calculus for realizing objects and morphism along a homological functor $F_*$ with enough injectives. In general we follow \cite{BlDG:pi-algebra} and \cite{GoHop:moduli}. See also \cite{Baues:combinatorial}. Since we are in a completely linear or stable situation the theory required to set up the obstruction calculus simplifies compared to the other settings. 
Nevertheless the simplifications in paragraph \ref{subsection:Moduli} compared to \cite{BlDG:pi-algebra} result from the use of truncated resolution model structures.
An obstruction calculus for realizing objects using only the triangulated structure is described, among other things, in \cite{BKS:Tate}.
We apologize in advance for using so much notation, but it seems unavoidable.

Our task was to look out for realizations in $\ho{\mc{M}}=\mc{T}$ of objects in \mc{A}. To motivate our next definition, let $X$ be an object in \mc{M} and let $X^\bullet\to r^0X$ be an $n$-\mc{G}-cofibrant approximation. We know:
     $$ \pi_s\stabhom{r^0X}{G}=\left\{
                               \begin{array}{cl}
                                   \stabhom{X}{G} &, \hbox{ if } s=0 \\ 
                                        0         &, \hbox{ else} 
                               \end{array} \right.  $$
With the spiral exact sequence we can calculate:
\begin{equation*} 
        \naturalpi{s}{X^\bullet}{G}=\left\{
                               \begin{array}{cl}
                                   \stabhom{X}{\Omega^s G} &,\hbox{ if } 0\sth s\sth n \\ 
                                        0         &,\hbox{ for } s\uber n 
                               \end{array} \right. 
\end{equation*}
And respectively:
     $$ \pi_s\stabhom{X^\bullet}{G}=\left\{
                               \begin{array}{cl}
                                   \stabhom{X}{G} &, \hbox{ if } s=0 \\ 
                                   \stabhom{X}{\Omega^{n+1}G} &,\hbox{ if } s=n+2 \\ 
                                        0         &,\hbox{ else} 
                               \end{array} \right. $$
All these sets of isomorphisms determine each other.
Of course, this is not the way we will encounter such spaces, since we are seeking for realization and not starting with them. Instead we will take these equations as the defining conditions of our successive realizations. 

\begin{Def} \label{potentielles n-Stück}
Let $A$ be an object in the abelian target category \mc{A}. We will call a Reedy cofibrant object $X^\bullet$ in $c\mc{M}$ a {\bf\boldmath potential $n$-stage for $A$\unboldmath} following \cite{BlDG:pi-algebra} and \cite{GoHop:moduli}, if there are natural isomorphisms
     $$ \naturalpi{s}{X^\bullet}{G}\cong\left\{
                               \begin{array}{cl}
                                   \Hom{\mc{A}}{A}{F_{*+s}G} &,\hbox{ if } 0\sth s\sth n \\ 
                                        0         &,\hbox{ for } s\uber n 
                               \end{array} \right. , $$
where we consider both sides as functors on \mc{G} as a subcategory of \mc{T}.
Note that this also makes sense for $n=\infty$. In this case an object satisfying these equations is simply called an {\bf\boldmath$\infty$\unboldmath-stage}. The reason is that by \ref{Unendlich-Stücke sind Realisierungen} it is not ``potential'' any more.
\end{Def}

\begin{Bem} 
If $X^\bullet$ is a potential $n$-stage for an object $A$ in \mc{A} then $\sk_nX^\bullet$ is a potential $\hrl n\!-\!1\hrr$-stage for $A$.
\end{Bem}

\begin{Bem} \label{Kohomotopie eines $n$-Stücks} 
Since \mc{G} was the class of $F$-injectives, the class $\kl F_*G|G\aus\mc{G}\kr$ is cogenerating the category \mc{A}, and we derive for a potential $n$-stage $X^\bullet$ from the previous properties and the spiral exact sequence the following equations:
     $$  \pi^sF_*X^\bullet\cong H^sNF_*X^\bullet
                          =\left\{
                               \begin{array}{cl}
                                   A     &, \hbox{ if } s=0 \\ 
                           A\hel n+1\her &, \hbox{ if } s=n+2 \\ 
                                   0     &,\hbox{ else }
                               \end{array} \right.  $$ 
The shift functor $\hel\frei\her$ is the internal shift from \ref{Einhängungen}.
\end{Bem}

Following the outlined philosophy we start the process of realizing an object $A$ in \mc{A} with a potential $0$-stage. 
Then we proceed by gluing on special objects to get to higher $n$-stages.
We define these layers in \ref{L(N,n)}. 
They have a certain representation property, see \ref{Darstellungseigenschaft von L(N,n)} and \ref{Kohomologie} will provide the obstruction groups, we are looking for.
To prove this property we have to consider algebraic analogues of these layers defined in \ref{K_A} and they should not be confused with each other. We also describe the moduli space of these different sorts of objects. In particular we will see, that they are connected, which ensures that the layers all look alike.

\begin{Def} \label{K_A}\label{K(N,n)}
Let $N$ be an object of \mc{A} and let $n\bth 0$.
We call an object $I^\bullet$ in $c\mc{A}$ an {\bf object of type \boldmath$K\hrl N,n\hrr$\unboldmath}, if the following conditions are satisfied:
    $$ \pi^sI^\bullet\cong\left\{
                   \begin{array}{cl}
                             N     &, \hbox{ if } s=n \\
                             0     &, \hbox{ else }
                   \end{array} 
                                        \right. $$
We denote $I^\bullet$ by $K(N,n)$. These objects are essentially unique by the following remark. For a quick introduction to moduli spaces we refer to appendix \ref{appendix:moduli}.
\end{Def}

\begin{Bem} \label{Modulraum von K(A,0)}
If $I^\bullet$ is an object of type $K\hrl A,0\hrr$, then there is a weak equivalence $r^0A=r^0\hrl\pi^0I^\bullet\hrr\to I^\bullet$. It follows that the moduli space is weakly equivalent to $B{\rm Aut}\hrl A\hrr$.

Objects of type $K\hrl N,n\hrr$ exist, for example is $\Omega^n_{\rm ext} r^0N$ or equivalently $\ol{W}\mbox{}^nr^0N$ such an object. The moduli space is given $B\text{Aut}\hrl N\hrr$, since the functor $\Omega_{\rm ext}$ induces an obvious equivalence $\mc{M}\hrl K\hrl N,n\hrr\hrr\to \mc{M}\hrl K\hrl N,n+1\hrr\hrr$ for $n\bth 0$. In particular this space is connected. 
\end{Bem}

\begin{Def} \label{L(N,n)}
Let $N$ be an object in \mc{A} and $n\bth 0$. We call an object $Y^\bullet$ in $c\mc{M}$ an {\bf object of type \boldmath$L\hrl N,n\hrr$\unboldmath}, if the following conditions are satisfied:
    $$ \naturalpi{s}{Y^\bullet}{G}=\left\{
                   \begin{array}{cl}
                       \Hom{\mc{A}}{N}{F_*G} &, \hbox{ if } s=n \\
                                 0           &, \hbox{ else }
                   \end{array} 
                                        \right. $$
Their existence and moduli space is described in \ref{Existenz von L(A,0)}.
We denote a generic object by $L(N,n)$. Do not confuse these objects with objects of type $K(N,n)$ in $c\mc{A}$, see \ref{K(N,n)}, the end of the remark \ref{Rechnung mit Spiralsequenz} and \ref{Konstruktion von phi}.
\end{Def}

\begin{Bem} \label{Existenz von L(A,0)}
Objects of type $L\hrl A,0\hrr$ exist: We choose an exact sequence
    $$ 0\to A\to I^0\toh{d} I^1 $$
with $I^0$ and $I^1$ injective. The map $d$ is induced by a map $E\hrl I^0\hrr\to E\hrl I^1\hrr$ in $\mc{T}=\ho{\mc{M}}$ between $F$-injective objects by \ref{Eilenberg-MacLane} that we will also call $d$. This $d$ again is represented by a map $d$ in \mc{M} if we choose the models for $E\hrl I^0\hrr$ and $E\hrl I^1\hrr$ to be fibrant and cofibrant. Now define a $1$-truncated cosimplicial object
\diagr{ E\hrl I^0\hrr \ar@<4pt>[r]\ar@<-4pt>[r] & E\hrl I^0\hrr\times E\hrl I^1\hrr \ar[l] }
with
\begin{equation*}
    d^0 = \left(  \begin{array}{c}                    
                   1  \\                                 
                   d \end{array} \right)  ,\
    d^1 = \left(  \begin{array}{c}                    
                   1  \\                                 
                   0 \end{array} \right)  ,\ \text{ and }
    s^0 = \left(  \begin{array}{cc}                    
                   1 & 0 \end{array} \right) . 
\end{equation*}
By applying a Reedy cofibrant approximation and left Kan extension we get a entire cosimplicial object, which is of type $L\hrl A,0\hrr$.

Objects of type $L\hrl N,n\hrr$ exist, since they can be given by setting: 
   $$ L\hrl N,n\hrr:=\Omega^n_{\rm ext} L\hrl N,0\hrr\ \text{ or }\ L\hrl N,n\hrr:=\ol{W}\mbox{}^n L\hrl N,0\hrr $$

We will compute the moduli space of objects of type $L\hrl A,0\hrr$ in \ref{Modulraum von L(A,0)}, it is given by ${\rm BAut}\hrl A\hrr$. 
Then the moduli space of objects of type $L\hrl N,n\hrr$ is given by $B\text{Aut}\hrl N\hrr$. This is proved by observing that it follows from \ref{Omega und W-quer} and \ref{noch mehr stabil} that $\Sigma_{\rm ext}$ and $\ol{W}$ induce mutually inverse homotopy equivalences of $\mc{M}\hrl L\hrl N,n\hrr\hrr$ and $\mc{M}\hrl L\hrl N,n+1\hrr\hrr$ for $n\bth 0$.
From lemma \ref{noch mehr stabil} we also get that $ \Sigma_{\rm ext} L\hrl N,n+1\hrr\cong \ L\hrl N,n\hrr$. 
\end{Bem}

\begin{Bem} \label{Rechnung mit Spiralsequenz}
By the spiral exact sequence we compute from \ref{L(N,n)}:
    $$ \pi_s\stabhom{L\hrl N,n\hrr}{G}=\left\{
                   \begin{array}{ll}
                       \Hom{\mc{A}}{N}{F_*G} &, \hbox{ if } s=n \\
                       \Hom{\mc{A}}{N}{F_{*+1}G} &, \hbox{ if } s=n+2 \\
                       \hspace{1cm}          0           &, \hbox{ else }
                   \end{array} 
                                        \right. $$
By the defining property of the $F$-injective objects in \mc{G} we get:
    $$ \pi^sF_*L\hrl N,n\hrr=\left\{
                   \begin{array}{cl}
                             N     &, \hbox{ if } s=n \\
                       N\hel 1\her &, \hbox{ if } s=n+2 \\
                             0     &, \hbox{ else }
                   \end{array} 
                                        \right. $$
Both sets of data are equivalent to the defining equations of an object of type $L\hrl N,n\hrr$ in Definition \ref{L(N,n)}. In particular it follows that $F_*L\hrl N,n\hrr$ is not an object of type $K\hrl N,n\hrr$.
We point out that for an object of type $L\hrl A,0\hrr$ the image $F_*L\hrl A,0\hrr$ is $1$-\mc{I}-equivalent to $r^0A$, as we can see with these equations.
\end{Bem}

\begin{Bem} \label{Konstruktion von phi}
Despite of the fact that $F_*L\hrl N,n\hrr$ is not of type $K\hrl N,n\hrr$ in $c\mc{A}$, there is a close connection explained in the following lemma. First we have to prepare ourselves. Let $n\bth 1$ and let $N$ be an object in \mc{A}. The isomorphism of $N$ and $\pi^nF_*L\hrl N,n\hrr$ defines a map $K\hrl N,n\hrr\to F_*L\hrl N,n\hrr$.
So by first applying $F_*$ and then pulling back along this arrow we obtain a map
\begin{equation} \label{Zurückziehen} 
    \phi_n\hrl Y^\bullet\hrr:{\rm map}\hrl L\hrl N,n\hrr, Y^\bullet\hrr \to \map\hrl K\hrl N,n\hrr, F_*Y^\bullet\hrr .
\end{equation}
Here we assume that $L\hrl N,n\hrr$ and $K\hrl N,n\hrr$ are Reedy cofibrant.
\end{Bem}

The next lemma is one of the central ingredients in the obstruction calculus as well as for the proof of \ref{Äquivalenz der 0-ten Schicht}.
\begin{lemma} \label{Darstellungseigenschaft von L(N,n)}
For $F$-fibrant $Y^\bullet$ in $c\mc{M}$ and Reedy cofibrant objects of type $L\emph{\hrl} N,n\emph{\hrr}$ and $K\emph{\hrl} N,n\emph{\hrr}$ the map $\phi_n\emph{\hrl} Y^\bullet\emph{\hrr}$ from \emph{\Ref{Zurückziehen}} is a natural weak equivalence.
\end{lemma}

\begin{beweis}
The proof is exactly parallel to the proof of \cite[Prop. 8.7.]{BlDG:pi-algebra}, although in our case linearity assures the result also for $n=0,1$. 
\end{beweis}

As to be expected, $0$-stages and $0$-layers will coincide. 
\begin{kor} \label{L(A,0) und 0-Stücke}
An object $X^\bullet$ is of type $L\emph{\hrl} A,0\emph{\hrr}$ if and only if it satisfies the following conditions:
\begin{punkt}
    \item There is an isomorphism $\pi^0F_*X^\bullet\cong A$ in \mc{A}.
    \item For every $Y^\bullet$ in $c\mc{M}$ the natural map
          $$ \emph{\stabhom{X^\bullet}{Y^\bullet}} \to \emph{\Hom{\mc{A}}{A}{\pi^0F_*Y^\bullet}} $$
is an isomorphism. 
\end{punkt}
\end{kor}

\begin{beweis}
Objects that satisfy (i) and (ii) are of type $L\hrl A,0\hrr$ because we can calculate:
\begin{align*}
    \naturalpi{s}{X^\bullet}{G}\cong\stabhom{X^\bullet}{\Omega^s_{\rm ext}r^0G}_F \cong \left\{ \begin{array}{cl}
                 \Hom{\mc{A}}{A}{F_*G} &, \text{ for } s=0 \\
                          0            &, \text{ else}
                 \end{array} \right.
\end{align*}
The other direction follows from \ref{Darstellungseigenschaft von L(N,n)}.
\end{beweis}

Finally we can determine the moduli space of all objects of type $L\hrl A,0\hrr$.
\begin{kor} \label{Modulraum von L(A,0)}
The moduli space of all objects of type $L\emph{\hrl} A,0\emph{\hrr}$ is connected and we have the following weak equivalence:
    $$ \mc{M}_F\emph{\hrl} L\emph{\hrl} A,0\emph{\hrr}\emph{\hrr} \simeq B {\rm Aut}\emph{\hrl} A\emph{\hrr} $$
\end{kor}

\begin{beweis}
The moduli space is connected:
Let $L\hrl A,0\hrr$ be some reference object and let $X^\bullet$ be another object of type $L\hrl A,0\hrr$. By pulling back $\id_A$ along the isomorphism of \ref{L(A,0) und 0-Stücke}(ii) we obtain a map $X^\bullet\to L\hrl A,0\hrr$ which induces an isomorphism on $\naturalpi{0}{\frei}{G}$ for every $G\aus\{F\text{-Inj}\}$. Both are potential $0$-stages, so this is the only group to check.

Now we will prove that the moduli space is weakly equivalent to the moduli space of objects of type $K\hrl A,0\hrr$. Then the result will follow from \ref{Modulraum von K(A,0)}.
By \ref{DK-Charakterisierung von Modulraeumen} there are canonical weak equivalences
    $$ \mc{M}_F\hrl L\hrl A,0\hrr\hrr \simeq B\text{haut}_F\hrl L\hrl A,0\hrr\hrr \hbox{ and } \mc{M}_{\mc{I}}\hrl K\hrl A,0\hrr\hrr \simeq B\text{haut}_{\mc{I}}\hrl K\hrl A,0\hrr\hrr. $$
It suffices to prove that ${\rm haut}_F\hrl L\hrl A,0\hrr\hrr\simeq {\rm haut}\hrl K\hrl A,0\hrr\hrr$, because both objects are fibrant grouplike simplicial monoids and $B$ preserves weak equivalences between fibrant simplicial sets. By \ref{Darstellungseigenschaft von L(N,n)} we have the following weak equivalences:
    $$ \map\hrl L\hrl A,0\hrr, L\hrl A,0\hrr\hrr\simeq \map\hrl K\hrl A,0\hrr, F_*L\hrl A,0\hrr\hrr\simeq \ell_0\Hom{\mc{A}}{A}{A} $$
Passing to appropriate components we see that ${\rm haut}_F\hrl L\hrl A,0\hrr\hrr\simeq \ell_0{\rm End}_{\mc{A}}\hrl A\hrr$ which finishes the proof.
Here $\ell_0\hrl\frei\hrr$ denotes the constant simplicial object.
\end{beweis}

\begin{Def} \label{Kohomologie}
Consider $K\hrl N,n\hrr$ in $c\mc{A}$ for $n\bth 0$. 
We assume that $K\hrl N,n\hrr$ is Reedy cofibrant.
Let $\Lambda^\bullet$ be an object in $c\mc{A}$. Then we define
    $$ \map\hrl K\hrl N,n\hrr, \wt{\Lambda}^\bullet\hrr =: \mc{H}^n\hrl \Lambda^\bullet,N\hrr, $$
where $\Lambda^\bullet\to\wt{\Lambda}^\bullet$ is a fibrant approximation, to be the {\bf\boldmath$n$\unboldmath-th cohomology space} of $\Lambda^\bullet$ with coefficients in $N$. We define the {\bf\boldmath$n$\unboldmath-th cohomology} of $\Lambda^\bullet$ by 
    $$ \pi_0\mc{H}^n\hrl \Lambda^\bullet,N\hrr  =: H^n\hrl \Lambda^\bullet,N\hrr.  $$
In the next lemma we will give an interpretation of these cohomology groups.
\end{Def}

\begin{Bem}
It follows for any $\Lambda^\bullet$ in $c\mc{A}$:
    $$ \Omega\mc{H}^n\hrl \Lambda^\bullet,N\hrr \simeq \mc{H}^{n-1}\hrl \Lambda^\bullet,N\hrr $$
\end{Bem}

\begin{lemma} \label{Kohomologie=Ext}
Let $\Lambda^\bullet$ in $c\mc{A}$ be \mc{I}-fibrant and $n$-\mc{I}-equivalent to $r^0\pi^0\Lambda^\bullet$.  
Then there is a natural isomorphism  
    $$ H^n\emph{\hrl} \Lambda^\bullet, N\emph{\hel} k\emph{\her}\emph{\hrr} \cong {\rm Ext}^{n,k}\emph{\hrl} N,\pi^0\Lambda^\bullet\emph{\hrr} $$
of abelian groups.  
\end{lemma}

\begin{beweis}
The canonical map $r^0\pi^0\Lambda^\bullet\to\Lambda^\bullet$ obtained by adjunction factors as the composition $r^0\pi^0\Lambda^\bullet\to\sk_{n+1}\Lambda^\bullet\to\Lambda^\bullet$ of $n$-\mc{I}-equivalences and we can approximate $\sk_{n+1}\Lambda^\bullet$ \mc{I}-fibrantly by $I^\bullet$ which yields an injective resolution of $\pi^0\Lambda^\bullet$ after normalization. Now $K\hrl N\hel k\her,n\hrr$ is an $\hrl n\!+\!1\hrr$-skeleton and we compute:
\begin{align*}
      H^n\hrl \Lambda^\bullet,N\hel k\her\hrr &= \pi_0\map\hrl K\hrl N\hel k\her,n\hrr, \Lambda^\bullet\hrr \cong \pi_0\map\hrl K\hrl N\hel k\her,n\hrr, I^\bullet\hrr \\
         &\cong {\rm Ext}^{n,k}\hrl N,\pi^0\Lambda^\bullet\hrr 
\end{align*}
\end{beweis}

\begin{Bem} \label{Bemerkung zu Ext}
Let $Y^\bullet$ be $F$-fibrant, such that $F_*Y^\bullet$ is $n$-\mc{I}-equivalent to $r^0\pi^0F_*Y^\bullet$.
Altogether lemma \ref{Darstellungseigenschaft von L(N,n)} and lemma \ref{Kohomologie=Ext} yield the following isomorphism of abelian groups:
    $$ \pi_0\map\hrl L\hrl A\hel k\her, n\hrr, Y^\bullet\hrr \cong \Ext{\mc{A}}{n,k}{A}{\pi^0F_*Y^\bullet} $$
Here we assume $L\hrl A\hel k\her, n\hrr$ and $Y^\bullet$ to be both Reedy cofibrant. This is functorial in $Y^\bullet$. It is not quite functorial in $A$, but for a morphism $A\to B$ after having chosen two objects $L\hrl A\hel k\her,n\hrr$ and $L\hrl B\hel k\her,n\hrr$ there is a uniquely determined homotopy class $L\hrl A\hel k\her,n\hrr\to L\hrl B\hel k\her,n\hrr$ inducing $A\to B$. The result tells us that an object $L\hrl N,n\hrr$ represents the cohomology functor $H^n\hrl F_*\hrl\frei\hrr,N\hrr$ in the homotopy category $\ho{c\mc{M}^F}$. Note that the isomorphism is in particular valid if $Y^\bullet$ is an $F$-fibrant $n$-stage.
\end{Bem}

We want to construct an obstruction calculus for lifting things from an interpolation category to the next one. In order to carry this out, we study the difference between potential $(n\!-\!1)$-stages and potential $n$-stages. 
In \ref{Pushout für ein n-Stück} and \ref{Konstruktion eines $n$-Stücks} we will prove the existence of certain homotopy pushout diagrams, where the difference between two stages is recognized as objects of type $L\hrl N,n\hrr$ for suitable $N$ and $n$. This construction can be viewed as a (potential) Postnikov cotower, compare \ref{ko-k-Invariante}.

The following two lemmas are an example for the simplifications we get for the stable case. The next lemma is the collapsed version of the so-called difference construction in \cite[8.4.]{BlDG:pi-algebra}.
\begin{lemma} \label{Kofaserkonstruktion in cM}
Let $n\bth 1$ and provide $c\mc{M}$ with the $F$-structure. Let $f:X^\bullet\to Y^\bullet$ be a map in $c\mc{M}$, which induces an isomorphism on $\pi^{0}F_*$ and whose homotopy cofiber $C^\bullet$ has the property that $\pi^sF_*C^\bullet=0$ for $0\sth s\sth n-1$. 
Let $P^\bullet$ be the homotopy fiber of $f$. Then there are isomorphisms $ P^\bullet\cong\Omega_{\rm ext}C^\bullet$ and $\Sigma_{\rm ext}P^\bullet\cong C^\bullet $ in $\ho{c\mc{M}^F}$ and $\sk_{n+2}P^\bullet$ is an object of type $L\emph{\hrl} \pi^nF_*C^\bullet, n+1\emph{\hrr}$.
\end{lemma}

\begin{beweis}
This follows directly from \ref{fast stabil} and \ref{noch mehr stabil} and the long exact $\pi_{*}^{\natural}$-sequence.
\end{beweis}

\begin{lemma} \label{Pushout für ein n-Stück}
Let $X^\bullet_n$ be a potential $n$-stage for $A$. Then $\sk_nX^\bullet_n=:X^\bullet_{n-1}$ is a potential $(n\!-\!1)$-stage for $A$, and there is a homotopy cofiber sequence in $c\mc{M}^F$:
     $$ L\emph{\hrl} A\emph{\hel} n\emph{\her},n+1\emph{\hrr} \to X^\bullet_{n-1} \to X^\bullet_n  $$
This sequence is also a homotopy fiber sequence in $c\mc{M}^F$.
\end{lemma}

\begin{beweis}
Call $C_n={\rm hocofib}\hrl X^\bullet_{n-1}\to X^\bullet_n\hrr$.
We know that $\pi_{s}^\natural$ of $C_n$ vanishes except in dimension $n$. 
Hence $\sk_{n+2}C_n$ is $F$-equivalent to $C_n$.
From \ref{Kofaserkonstruktion in cM} we see, that $\Sigma_{\rm ext}\Omega_{\rm ext}C_n\simeq C_n\simeq\Omega_{\rm ext}\Sigma_{\rm ext}C_n$ and that $\Omega_{\rm ext}C_n$ is an object of type $L\hrl A\hel n\her,n+1\hrr$. 
We also see that the sequence is a homotopy cofiber sequence as well as a homotopy fiber sequence.
\end{beweis}

\begin{lemma} \label{Konstruktion eines $n$-Stücks}
Let there be given a homotopy cofiber sequence in $c\mc{M}^F$:
\diagr{ L\emph{\hrl} A\emph{\hel} n\emph{\her},n+1\emph{\hrr} \ar[r]^-{w_{n}} & X^\bullet_{n-1} \ar[r] & X^\bullet_n  }
Let $X^\bullet_{n-1}$ be a potential $(n\!-\!1)$-stage for $A$.  
A Reedy cofibrant approximation to $X^\bullet_n$ is a potential $n$-stage for $A$ if and only if the map $w_{n}$ induces an isomorphism $ A\emph{\hel} n\emph{\her}\cong \pi^{n+1} F_*X^\bullet_{n-1}$.
\end{lemma}

\begin{beweis}
It follows from \ref{lange exakte Homologiesequenz} that there is an exact sequence 
    $$ 0\to \pi^nF_*X^\bullet_n\to\pi^{n+1}F_*L\hrl A\hel n\her,n+1\hrr\stackrel{\cong}{\longrightarrow}\pi^{n+1}F_*X^\bullet_{n-1}\to\pi^{n+1}F_*X^\bullet_n\to 0 $$
and an isomorphism $\pi^{n+2}F_*X^\bullet_n\cong\pi^{n+3}F_*L\hrl A\hel n\her,n+1\hrr\cong A\hel n+1\her$. All other groups of the form $\pi^sF_*X^\bullet_n$ for $s\uber 0$ vanish, hence $X^\bullet_n$ has the right homotopy groups for a potential $n$-stage for $A$, we just need to approximate it Reedy cofibrantly.
\end{beweis}

\begin{Def} \label{ko-k-Invariante}
We will call a map $w_n$ as in \ref{Konstruktion eines $n$-Stücks} an {\bf\boldmath $n$-th attaching map}, i.e. if it is an $F$-cofibration of the form $L\hrl A\hel n\her, n+1\hrr\to X^\bullet_{n-1}$ between Reedy cofibrant objects, whose target is a potential $n$-stage. The induced homotopy class will be called {\bf\boldmath $n$-th co-k-invariant}.
The concept is dual to that of $k$-invariants of a Post\-ni\-kov-tower. 

If we need to refer to an attaching map stemming from a specified potential $n$-stage $X_n^\bullet$ as in \ref{Pushout für ein n-Stück}, we will write $w^{X^\bullet_n}$ instead of $w_n$.
To an $\infty$-stage we can associate attaching maps for each of its $n$-stages. By abuse of notation we will denote them by $w^{X^\bullet_n}$ without specifying an actual $n$-stage.
\end{Def}

Now we start to describe the obstruction against the existence of realizations of objects in $IP_{n-1}\hrl F\hrr$.

\begin{Def} 
Let $X^\bullet_{n-1}$ be a potential $(n-1)$-stage for an object $A$. We call an object $X^\bullet$ a {\bf\boldmath potential $n$-stage over $X^\bullet_{n-1}$\unboldmath} if $X^\bullet$ is a potential $n$-stage and $\sk_nX^\bullet$ is $F$-equivalent to $X^\bullet_{n-1}$. This is equivalent to $X^\bullet_{n-1}$ being $\hrl n-1\hrr$-$F$-equivalent to $\sk_nX^\bullet$.
\end{Def}

The obstruction against the existence of an $n$-stage over a given $(n-1)$-stage is the existence of an attaching map $w_n$ like in \ref{Konstruktion eines $n$-Stücks}. We are now going to reformulate this in algebraic terms. We already know from remark \ref{Kohomotopie eines $n$-Stücks}  that for an $(n\!-\!1)$-stage $X_{n-1}^\bullet$ its image $F_*X_{n-1}^\bullet$ has the same cohomology groups as an object of type $K\hrl A,0\hrr\oplus K( A\hel n\her,n+1)$. 
Without loss of generality we assume $X^\bullet_{n-1}$ to be $F$-fibrant.
Hence we know that such an attaching map $w_n$ exists if and only if we are able to construct a map
    $$ \omega_{n}:K\hrl A\hel n\her,n+1\hrr \to F_*X_{n-1}^\bullet $$
inducing an isomorphism on $\pi^{n+1}F_*\hrl\frei\hrr$, because by the representing property \ref{Darstellungseigenschaft von L(N,n)} it follows that we were then able to choose $w^{X^\bullet_n}$ such that 
   $$\pi_0\hel\phi\hrl X^\bullet_{n-1}\hrr\hrl w_n\hrr\her= \pi_0\hel\omega_{n}\her.$$
From \ref{Kohomotopie eines $n$-Stücks} we have the homotopy cofiber sequence
\begin{equation*} 
     K\hrl A\hel n\her,n+2\hrr \to \sk_{n+1}F_*X^\bullet_{n-1} \to F_*X^\bullet_{n-1} 
\end{equation*}
and we can consider the following diagram:
\diag{ & \sk_1F_*X_{n-1}^\bullet \ar[r]^-{\simeq}\ar[d]^{\cong} & K\hrl A,0\hrr\simeq r^0A \\
        K\hrl A\hel n\her,n+2\hrr \ar[r]^-{\beta_{n}}\ar[d] \ar@{}[dr]|->>>{\phantom{xx}{\rm ho}-\pushout} & \sk_{n+1}F_*X^\bullet_{n-1} \ar[d] & \\
        \ast  \ar[r]\ar[d] & \sk_{n+2}F_*X^\bullet_{n-1} \ar[r]^-{\cong} & F_*X_{n-1}^\bullet \\
        K\hrl A\hel n\her,n+1\hrr \ar@{.>}@/_13pt/[rru]_-{\omega_{n}} &&  }{b-Invariante} 
Observe also that we have isomorphisms
    $$ H^{n+2}\hrl \sk_{n+1}F_*X_{n-1}^\bullet, A\hel n\her\hrr \stackrel{\cong}{\longrightarrow} H^{n+2}\hrl r^0A, A\hel n\her\hrr=\Ext{\mc{A}}{n+2,n}{A}{A} $$
of abelian groups by lemma \ref{Kohomologie=Ext} or remark \ref{Bemerkung zu Ext}.
 
\begin{Def} \label{Definition b-Invariante}
The homotopy class $b_{n}$ of the map $\beta_{n}$ in 
    $$\pi_0\map\hrl K\hrl A\hel n\her,n+2\hrr,r^0A\hrr= H^{n+2}\hrl r^0A, A\hel n\her\hrr\cong\Ext{\mc{A}}{n+2,n}{A}{A} $$ 
will be called the {\bf obstruction class} of the potential $\hrl n-1\hrr$-stage $X^\bullet_{n-1}$. 
\end{Def}

\begin{lemma} \label{b-Invariante nullhomotop}
In \emph{\Ref{b-Invariante}} the map $\omega_{n}$ inducing an isomorphism on $\pi^{n+1}F_*\emph{\hrl}\frei\emph{\hrr}$ exists if and only if $\beta_{n}$ is nullhomotopic.
\end{lemma}

\begin{beweis}
Obvious.
\end{beweis}

\begin{satz} \label{Ext^(n+2,n) und Objekte}
Let $n\bth 1$ and $A$ be an object of \mc{A}. Let $X^\bullet_{n-1}$ be a potential $(n-1)$-stage of $A$. 
There exists a potential $n$-stage $X^\bullet_n$ over $X^\bullet_{n-1}$ if and only if the co-$k$-invariant $b_{n}$ from definition \emph{\ref{Definition b-Invariante}} in ${\rm Ext}^{n+2,n}\emph{\hrl} A,A\emph{\hrr}$ vanishes.
\end{satz}

\begin{beweis}
From \ref{b-Invariante nullhomotop} we know, that an attaching map for $X^\bullet_n$ exists if and only if $b_n=0$.
\end{beweis}

Now we are concerned with telling apart different realizations.

\begin{Def} \label{Differenzklasse für Objekte}
Let $X^\bullet_n$ and $Y^\bullet_n$ be potential $n$-stages for an object $A$ with $\sk_n X^\bullet_n \simeq X^\bullet_{n-1} \simeq \sk_n Y^\bullet_n$. The homotopy fiber of the canonical maps from $X^\bullet_{n-1}$ to $X^\bullet_n$ and $Y^\bullet_n$ is $L\hrl A\hel n\her,n+1\hrr$ by \ref{Pushout für ein n-Stück}. We obtain two attaching maps
    $$ w^{X^\bullet_n} \text{ and } w^{Y^\bullet_n}:L\hrl A\hel n\her,n+1\hrr \to X^\bullet_{n-1}, $$
The {\bf difference class} of the objects $X^\bullet_n$ and $Y^\bullet_n$ is defined to be the class
    $$ \delta\hrl X^\bullet_n, Y^\bullet_n\hrr:= \pi_0\hrl w^{X^\bullet_n}\hrr - \pi_0\hrl w^{Y^\bullet_n}\hrr \aus \pi_0\mc{H}^{n+1}\hrl F_*X^\bullet_{n-1}, A\hel n\her\hrr\cong {\rm Ext}^{n+1,n}\hrl A,A\hrr . $$
\end{Def}

\begin{Bem} 
The proof of the next theorem shows that this defines an action of $\Ext{\mc{A}}{n+1,n}{A}{A}$ on the class of $F$-equivalence classes of potential $n$-stages over a given potential $\hrl n\!-\!1\hrr$-stage. It is obviously transitive. This proves first of all that there is just a set of such equivalence classes or, what is the same, of realizations in $IP_n\hrl F\hrr$ of a given object in $IP_{n-1}\hrl F\hrr$. 
\end{Bem}

\begin{satz} \label{Ext^(n+1,n) und Objekte}
Let $n\bth 1$.
There is an action of ${\rm Ext}^{n+1,n}\emph{\hrl}A,A\emph{\hrr}$ on the set of $F$-equivalence classes of potential $n$-stages of $A$ over a given potential $\emph{\hrl}n\!-\!1\emph{\hrr}$-stages, if it is non-empty, which is transitive and free. 
\end{satz}

\begin{beweis}
Let $X^\bullet_n$ be a potential $n$-stage with $X^\bullet_{n-1}:=\sk_n X^\bullet_n$ and take a class $\kappa\aus\Ext{\mc{A}}{n+1,n}{A}{A}$. We want to construct a potential $n$-stage $Y^\bullet_n$ over $X^\bullet_{n-1}$, such that $\kappa=\delta\hrl X^\bullet_n,Y^\bullet_n\hrr$.
Consider the map $\gamma$ given by the following composition 
    $$ \xymatrix{ K\hrl A\hel n\her,n+1\hrr\ar[r]^-{\kappa} & K\hrl A,0\hrr \ar[r]^-{\text{incl.}} &  K\hrl A,0\hrr\oplus K\hrl A\hel n\her, n+1\hrr \cong F_*X^\bullet_{n-1} },$$  
and let $c\co L\hrl A\hel n\her, n+1\hrr \to X^\bullet_{n-1}$ be a realization of $\gamma$ existing by \ref{Darstellungseigenschaft von L(N,n)}.
Take a map 
    $$\omega_{n}\co K\hrl A\hel n\her,n+1\hrr\to F_*X^\bullet_{n-1}$$
from \Ref{b-Invariante} representing the homotopy class of 
    $$w^{X^\bullet_{n}}\co L\hrl A\hel n\her,n+1\hrr\to X^\bullet_{n-1}$$
associated to $X^\bullet_n$ by \ref{Pushout für ein n-Stück} and add it to $c$. 
The resulting map $\omega_{n}+\gamma$ will still induce an isomorphism on $\pi^{n+1}$, since $\gamma$ itself induces the zero map on $\pi^{n+1}$.
Thus the cofiber $Y^\bullet_n$ of the corresponding map $w^{X^\bullet_n}+c\co L\hrl A\hel n\her,n+1\hrr\to X^\bullet_{n-1}$ is a potential $n$-stage over $X^\bullet_{n-1}$ by \ref{Konstruktion eines $n$-Stücks}, which realizes the given difference class, hence $\kappa=\delta\hrl X^\bullet_n,Y^\bullet_n\hrr$. 
This process is obviously additive in $[\kappa]$, therefore we have a group action. It is also clear that $X^\bullet_n\cong Y^\bullet_n$ in $IP_n\hrl F\hrr$ if and only if $\kappa=0$.
\end{beweis}



Now we are going to describe the obstruction for lifting maps. 
\begin{Def} \label{Liftung}
Let $n\bth 1$. Let $X^\bullet$ and $Y^\bullet$ be objects in $IP_n\hrl F\hrr$ and let $\varphi:\sigma_nX^\bullet\to\sigma_nY^\bullet$ be a map in $IP_{n-1}\hrl F\hrr$. We say that {\bf\boldmath $\varphi$ lifts \unboldmath} if there is a map $\Phi:X^\bullet\to Y^\bullet$ such that $\sigma_n\Phi=\varphi$.
In this case we call $\Phi$ a {\bf lifting} of $\varphi$.
\end{Def}

\begin{Bem} \label{Erklärungen zu Liftungen}
By definition every object $W^\bullet$ in $IP_n\hrl F\hrr$ can be approximated by a potential $n$-stage, which corresponds to $n$-$F$-cofibrant approximation. On the other side it can be approximated $F$-fibrantly such that $F_*W^\bullet$ is $\hrl n+1\hrr$-\mc{I}-equivalent to $r^0\pi^0F_*W^\bullet$.

Between a potential $n$-stage $X^\bullet$ and an $F$-fibrant $Y^\bullet$ where $F_*Y^\bullet$ is $\hrl n+1\hrr$-\mc{I}-equivalent to $\pi^0F_*Y^\bullet$ every morphism in $IP_n\hrl F\hrr$ can be represented by a map $f:X^\bullet\to Y^\bullet$ in $c\mc{M}$. 

Assume that we are given a morphism from $X^\bullet$ to $Y^\bullet$ in $IP_{n-1}\hrl F\hrr$, then this can be represented by a map $f:\sk_nX^\bullet\to Y^\bullet$. Now $f$ lifts if and only if there is a map $\wt{f}:X^\bullet\to Y^\bullet$ such that 
   $$ \xymatrix@=13pt{ \sk_nX^\bullet \ar[r]& X^\bullet \ar[r] & Y^\bullet } $$
is homotopic to $f$ in $c\mc{M}^F$. 
\end{Bem}

\begin{satz} \label{Ext^(n+1,n) und Morphismen}
A morphism $\sigma_nX^\bullet\to \sigma_nY^\bullet$ in $IP_{n-1}\emph{\hrl}F\emph{\hrr}$ lifts to a morphism $X^\bullet\to Y^\bullet$ in $IP_{n}\emph{\hrl}F\emph{\hrr}$ if and only if ${\rm ob}_n(f)\aus\emph{\Ext{\mc{A}}{n+1,n}{\pi^0F_*X^\bullet}{\pi^0F_*Y^\bullet}}$ defined in \emph{\Ref{ob}} vanishes.
\end{satz}

\begin{beweis}
We assume without loss of generality that $X^\bullet$ is a potential $n$-stage for an object $A$ and that $Y^\bullet$ is $F$-fibrant such that $F_*Y^\bullet$ is $\hrl n+1\hrr$-\mc{I}-equivalent to $r^0B$ in $c\mc{A}$. We can achieve this by approximations in the $n$-$F$-structure. 
Also without loss of generality we can replace to homotopy cofiber sequence  
\diagr{ L\hrl A\hel n\her,n+1\hrr \ar[r]^-{w^{X^\bullet_n}} & X^\bullet_{n-1} \ar[r] & X^\bullet_n }
in $c\mc{M}^F$ of \ref{Konstruktion eines $n$-Stücks} by an actual cofiber sequence using factorizations in the $F$-structure.   
This means that we have constructed the following solid arrow diagram
\diag{ L\hrl A\hel n\her,n+1\hrr \ar[r]^-{w^{X^\bullet_n}}\ar[d]\ar@{}[dr]|->>>{\pushout} & \sk_nX^\bullet \ar[r]^-f\ar[d] & Y^\bullet \\
       PL \ar[r] \ar@{.>}[urr] & X^\bullet \ar@{.>}[ur] &  }{Hindernis für f}
where $PL\stackrel{F}{\simeq}\ast$ is a path object in the $F$-structure for $L\hrl A\hel n\her,n+1\hrr$. We conclude that the existence of the dotted liftings in diagram \Ref{Hindernis für f} are equivalent to each other. 
By \ref{Darstellungseigenschaft von L(N,n)} we deduce that an extension of $f$ to $X^\bullet$ exists if and only if the map
   $$\xymatrix@=20pt{ K\hrl A\hel n\her, n+1\hrr \ar[rr]^-{\phi(fw^{X^\bullet_n})} && F_*Y^\bullet } $$
is null homotopic, where $\phi$ is the map from \Ref{Zurückziehen}. $\phi(fw^{X^\bullet_n})$ defines an obstruction element 
\begin{equation} \label{ob}
      {\rm ob}_n(f):=[\phi(fw^{X^\bullet_n})]\aus H^{n+1}\hrl F_*Y^\bullet, A\hel n\her\hrr =\Ext{\mc{A}}{n+1,n}{\pi^0F_*X^\bullet}{\pi^0F_*Y^\bullet}
\end{equation}
by \ref{Kohomologie=Ext}. Recall that $F_*Y^\bullet$ is $\hrl n\!+\!1\hrr$-\mc{I}-equivalent to $r^0\pi^0F_*Y^\bullet=r^0B$. So remark \ref{Bemerkung zu Ext} applies.
\end{beweis}

Before we proceed to the next theorem we have to reformulate the obstruction defined in \Ref{ob}. Here we use the ``almost stability'' of $\ho{c\mc{M}^F}$ that is displayed in \ref{fast stabil}, \ref{noch mehr stabil} and \ref{Pushout für ein n-Stück}. 
In the situation of \Ref{Hindernis für f} we want to achieve, that $f$ factors over a potential $\hrl n-1\hrr$-stage of $Y^\bullet$. Let
   $$ \xymatrix{ \sk_nY^\bullet \ar[r]^{\upsilon} & \wt{Y}^\bullet_{n-1}\ar[r]^{\wt{\upsilon}} & Y^\bullet} $$
be a factorization of the canonical inclusion map into an $\hrl n\!-\!1\hrr$-$F$-cofibration $\upsilon$ followed by an $\hrl n\!-\!1\hrr$-$F$-trivial fibration $\wt{\upsilon}$. The map $\upsilon$ will necessarily be an $F$-trivial cofibration. Obviously, $\wt{Y}^\bullet_{n-1}$ is a potential $\hrl n\!-\!1\hrr$-stage and the co-$k$-invariant of $Y^\bullet$ is represented by $\wt{\upsilon}\circ w^{\wt{Y}^\bullet_n}$. 
Then the map $f:\sk_nX^\bullet\to Y^\bullet$ factors over a map $\wt{f}:\sk_nX^\bullet\to\wt{Y}^\bullet_{n-1}$ whose homotopy class is uniquely determined. 
The attaching map $w^{Y^\bullet_n}$ prolongs to an attaching map
\diagr{ L\hrl B\hel n\her,n+1\hrr \ar[r]_-{w^{Y^\bullet_n}} & \sk_nY^\bullet \ar[r]^{\upsilon} & \wt{Y}^\bullet_{n-1} ,}
which we will denote by $w^{\wt{Y}^\bullet_n}$. Since $\upsilon$ is an $F$-equivalence, it induces the same class in $\Ext{}{n+1,n}{B}{B}$ as $w^{Y^\bullet_n}$. Consider the following diagram:
\diag{ L\hrl A\hel n\her, n+1\hrr \ar[r]^-{w^{X^\bullet_n}} \ar[d]_{\ol{f}} & \sk_nX^\bullet \ar[d]^{\wt{f}} \\
       L\hrl B\hel n\her, n+1\hrr \ar[r]_-{w^{\wt{Y}^\bullet_n}} & \wt{Y}^\bullet_{n-1}  }{anderes ob}
Here $\ol{f}$ is induced by $f$ in the following sense: It represents the uniquely determined homotopy class that induces the map $\pi^0F_*\hrl f\hrr:A\to B$, compare \ref{Bemerkung zu Ext}.

We observe that, if we are given a diagram like \Ref{Hindernis für f}, we get a diagram \Ref{anderes ob} and we have the following equation
\begin{equation}\label{Gleichheit der Homotopieklassen}
      \text{ob}_n\hrl f\hrr=[\phi(fw^{X^\bullet_n}]=\hrl\wt{\upsilon}\hrr_*\left( [w^{\wt{Y}^\bullet_n}\ol{f}]-[\wt{f}w^{X^\bullet_n}]\right) \ \aus\ \Ext{\mc{A}}{n+1,n}{A}{B} ,
\end{equation}
since $\xymatrix@1{L\hrl B\hel n\her, n+1\hrr \ar[r]^-{w^{\wt{Y}^\bullet_n}} & \wt{Y}^\bullet_{n-1}\ar[r]^{\wt{\upsilon}} & Y^\bullet}$ is a homotopy cofiber sequence.

\begin{lemma} \label{umformuliertes ob}
The obstruction ${\rm ob}_n(f)$ from \emph{\Ref{ob}} vanishes if and only if the diagram \emph{\Ref{anderes ob}} commutes in \ho{c\mc{M}^F}.
\end{lemma}

\begin{beweis}
If the square commutes up to homotopy we can strictify it by changing $\ol{f}$ and $\wt{f}$ within their homotopy class. Then we can apply the pushout functor and obtain a map $X^\bullet\to Y^\bullet$. We can turn this process around if we remember that the homotopy cofiber sequence in \ref{Pushout für ein n-Stück} is also a homotopy fiber sequence. So if a lifting exists, which is equivalent to ${\rm ob}_n(f)=0$, then this diagram commutes.
\end{beweis}

\begin{satz} \label{Derivation}
Let $X^\bullet$ and $Y^\bullet$ be objects in $IP_n\emph{\hrl}F\emph{\hrr}$ such that $\pi^0F_*X^\bullet=A$ and $\pi^0F_*Y^\bullet=B$. Then the map
    $$ {\rm ob}_n\emph{\hrl}\frei\emph{\hrr}\co \emph{\Hom{IP_{n-1}(F)}{\sigma_n X^\bullet}{\sigma_n Y^\bullet}}\to \emph{\Ext{\mc{A}}{n+1,n}{A}{B}} $$
is a homomorphism satisfying property {\rm (iii)} of \emph{\ref{exakte Sequenz von Kategorien}}. 
\end{satz}

\begin{beweis}
There is a map of sets $\Hom{IP_{n-1}(F)}{\sigma_n X^\bullet}{\sigma_n Y^\bullet}\to \Ext{\mc{A}}{n+1,n}{A}{B}$, where we map an $f$ like in \Ref{Hindernis für f} to $\text{ob}_n(f)$ as defined in \Ref{ob}. This is well defined by \ref{Darstellungseigenschaft von L(N,n)}. We easily see that it is a homomorphism of abelian groups when we put the fold map $Y^\bullet\oplus Y^\bullet\to Y^\bullet$ in the place of $Y^\bullet$ in diagram \ref{Hindernis für f}.

Property (iii) of \ref{exakte Sequenz von Kategorien} follows immediately by \ref{umformuliertes ob} by considering two squares like \Ref{anderes ob} for $[f]:\sigma_nX^\bullet\to \sigma_nY^\bullet$ and $[g]:\sigma_n Y^\bullet\to\sigma_n Z^\bullet$ respectively.
\diagr{ L\hrl A\hel n\her, n+1\hrr \ar[r]^-{w^{X^\bullet_n}} \ar[d]_{\ol{f}} & \sk_nX^\bullet \ar[d]^{\wt{f}} \\
       L\hrl B\hel n\her, n+1\hrr \ar[r]_-{w^{Y^\bullet_n}} \ar[d]_{\ol{g}} & \wt{Y}^\bullet_{n-1} \ar[d]^{\wt{g}} \\
       L\hrl C\hel n\her, n+1\hrr \ar[r]_-{w^{Z^\bullet_n}} & \wt{Z}^\bullet_{n-1}  }
Now the statement follows directly from \ref{Gleichheit der Homotopieklassen} and \ref{anderes ob}. The homotopy classes involving the term $w^{Y^\bullet_n}$ cancel out.
\end{beweis}

\begin{Bem} \label{Adams-Differential}
The homomorphism $\rm{ob}_n$ is identified by \ref{ASS} and \ref{Adams-Filtrierung} with the differential $d_n\co E_n^{0,0}\hrl X,Y\hrr\to\Ext{\mc{A}}{n+1,n}{FX}{FY}$ in the Adams spectral sequence for $F$. Here $E_n^{0,0}\hrl X,Y\hrr$ is the intersection of the kernels of the previous differentials. In particular we get from \Ref{Gleichheit der Homotopieklassen} the following formula
    $$ d_nf = w_{n-1}^Yf-fw_{n-1}^X $$
for $f\aus\Hom{\mc{A}}{FX}{FY}$, which is reminiscent of \cite[Prop. 8.10.]{Bou:klocal}. Here $w_{n-1}^X$ and $w_{n-1}^Y$ are the co-$k$-invariants of $X$ and $Y$.
\end{Bem}

From diagram \Ref{Hindernis für f} we are now going to derive the obstruction against the uniqueness of the realization of $f$. 

\begin{Def} \label{Operation von Ext^nn}
Let $f:X^\bullet\to Y^\bullet$ be a map of potential $n$-stages for $A$ and $B$ respectively and for $n\bth 1$. Let $\alpha\aus\Ext{\mc{A}}{n,n}{A}{B}$. 
We define a new map $\alpha+f:X^\bullet\to Y^\bullet$ by the following diagram:
\begin{equation} \label{DIAGRAMM}
\xy
   \xymatrix"*"@=21pt{ \ast \ar[dd] \ar[rr] & & L\hrl A\hel n\her,n\hrr \ar@{.>}@/^13pt/[ddrr]^-{\alpha} & & \\ 
                                & & & &              \\
                       X^\bullet \ar[rrrr]_-{\alpha+f}  & & & & Y^\bullet  }
   \POS(11,10)
   \xymatrix@=17pt{ L\hrl A\hel n\her,n+1\hrr  \ar[rr] \ar'[d] [dd] \ar["*"]   & & \ast \ar[dd] \ar["*"] \ar@{.>}@/^13pt/["*"ddrr]^-{f''}  \\ 
                                                & &                    \\
              \sk_n X^\bullet \ar[rr] \ar["*"] \ar@{.>}["*"rrrr]_-<<<<<<{f'} & & X^\bullet \ar["*"rr]_-<{f} }    \endxy
\end{equation}
This diagram commutes up to homotopy, it can be strictified by choosing appropriate replacements for $*$ by external cone objects.
The top square, the square on the left and the square in the back part of \Ref{DIAGRAMM} are homotopy pushout squares. 
The datum of a map $f$ is equivalent to giving maps $f'$ and $f''$ making the obvious square (homotopy) commutative. 
Prescribing the homotopy class of $\alpha$ is equivalent to the existence of a map $\alpha+f$ whose homotopy class is, like that of $f$, a lifting of the homotopy class of $f'$ in the sense of \ref{Liftung}. Of course, the homotopy class of $\alpha+f$ is uniquely determined by the homotopy class of $\alpha$ (and of $f$).
This defines a group action of $\Ext{\mc{A}}{n,n}{A}{B}$ on the liftings of a map $\sk_nX^\bullet\to Y^\bullet$ in $IP_n\hrl F\hrr$, if they exist. 
We will denote this action by $\gamma_n$, so that $\alpha+f=\gamma_n(\alpha,f)$.
\end{Def}

\begin{satz} \label{Ext^(n,n) und Morphismen}
The map $\gamma_n$ is an action of $\emph{\Ext{\mc{A}}{n,n}{A}{B}}$ on $\emph{\Hom{IP_n(F)}{X^\bullet}{Y^\bullet}}$. Two morphisms agree on the $n$-skeleton if and only if they are in the same orbit. The restriction of the action to the set of realizations in $IP_n\emph{\hrl} F\emph{\hrr}$ of a given morphism in $IP_{n-1}\emph{\hrl} F\emph{\hrr}$ is free and transitive. The restricted action satisfies the linear distributivity law from \emph{\ref{exakte Sequenz von Kategorien}}.  
\end{satz}

\begin{beweis}
All facts are straightforward.
To prove that the linear distributivity law holds we observe that the data to construct $\hrl\beta+g\hrr\hrl\alpha+f\hrr$ is contained in the map
\diagr{ L\hrl A\hel n\her,n\hrr \ar[r]^{\rm diag} & 
        *{\begin{array}{c} L\hrl A\hel n\her,n\hrr \\ \oplus \\ L\hrl A\hel n\her,n\hrr  \end{array}} \ar@/^20pt/[rr]^-{\alpha} \ar@/_15pt/[r]_-{\ol{f}} & 
        L\hrl B\hel n\her,n\hrr \ar@/_15pt/[rr]_-{\beta}
        & Y^\bullet \ar[r]^g & Z^\bullet ,}
where the homotopy class of $\ol{f}$ is induced by $f:X^\bullet\to Y^\bullet$. 
The induced map $X^\bullet\to Z^\bullet$ is $g_*\alpha+f^*\beta+gf$.
\end{beweis}

\subsection{The tower of interpolation categories}
\label{subsection:Turm von Interpolationskategorien}

We will now plug the $n$-$F$-struc\-tures from theorem \ref{(abgeschnittene) F-Strukturen} into the tower of truncated homotopy categories in \ref{Turm von abgeschnittenen Homotopiekategorien}. 
Suitable subcategories then supply a tower of interpolation categories for $F$. There are two ways to describe this tower. The first one is closer to the philosophy in \cite{BlDG:pi-algebra}.

\begin{Def} \label{Interpolationskategorien}
Let $n\bth 0$. Let $IM_n\hrl F\hrr$ be the full subcategory of $c\mc{M}$ that consists of those objects $X^\bullet$, which are $F$-equivalent to a potential $n$-stage for some $A$ in \mc{A} defined in \ref{potentielles n-Stück}.  
We call this category the {\bf\boldmath $n$-th interpolation model of $F$\unboldmath}.

Let $IP_n(F)$ be the image of $IM_n(F)$ in $\ho{c\mc{M}^{F}}$.
We call this category the {\bf\boldmath $n$-th interpolation category of the functor $F$\unboldmath}.
\end{Def}

\begin{Bem} \label{alternativ}
The second one is closer to our idea of truncated objects, where only the front end of the objects count and where we do not care for the higher degrees. Consider the full subcategory of $c\mc{M}$ of objects, which are $n$-$F$-equivalent to some potential $n$-stage. Then $IP_n(F)$ is equivalent to the image of this subcategory in $\ho{c\mc{M}^{n-F}}$.
The notion of isomorphism in $\ho{c\mc{M}^{0-F}}$ is rather coarse, hence a lot of objects become identified. As $n$ grows, fewer and fewer objects qualify for $IP_n\hrl F\hrr$, while the equivalences get finer and finer. 

Note that whether $X^\bullet$ belongs to some interpolation category or not is detected by $FX^\bullet$, see remark \ref{Kohomotopie eines $n$-Stücks}.
\end{Bem}

\begin{Bem} \label{F-Turm} 
The functors $\sigma_n$ from our tower of truncated homotopy categories restrict to our interpolation categories $IP_n\hrl F\hrr$. 
Also $\theta_nX$ for some $X$ in \mc{T} lands in $IP_n\hrl F\hrr$.
We will continue to denote these functors by $\theta_n$ and $\sigma_n$.

There is also an additional functor $\pi^0F_*\cong H^0NF_*:IP_n\hrl F\hrr\to\mc{A}$, which is derived from the functor $c\mc{M}\to\mc{A}$ given by
    $$ X^\bullet \mapsto \pi^0F_*X^\bullet .$$
We arrive at the following tower of interpolation categories:
\diag{  & \mc{T} \ar[dl]_{\theta_{n+1}}\ar[d]^{\theta_n}\ar[dr]^{...}\ar[rrr]^F &&& \mc{A} \\
         IP_{n+1}\hrl F\hrr \ar[r]_{\sigma_n} & IP_{n}\hrl F\hrr \ar[r] & \hdots  \ar[r] & IP_1\hrl F\hrr \ar[r]_{\sigma_0} & IP_0\hrl F\hrr \ar[u]_{H^0NF_*=\pi^0F_*}  }{Interpolationsturm}
Again this diagram $2$-commutes in the $2$-category of categories. It is worth to emphasize that the restricted functors $\sigma_n$ and $\theta_n$ here in general do not possess left adjoints contrary to the ones in \ref{Turm von abgeschnittenen Homotopiekategorien}. The left adjoints there do not take values in some interpolation category.
\end{Bem}

\subsection{Spectral sequences}
\label{subsection:SS}

Before we start the discussion about spectral sequences, let us explain what we want to find out. We would like to state that if we have found an object $X^\bullet$ which looks like a fibrant replacement of a constant object, thus $\pi^sFX^\bullet=0$ for all $s\uber 0$, then $\Tot X^\bullet$ is an actual realization of $\pi^0FX^\bullet$.
We will use this in \ref{die Abbildung Real_8(A) -> Real(A)} and \ref{Unendlich-Stücke sind Realisierungen}. 
Unfortunately there are problems with the spectral sequences in the cosimplicial case which do not occur in a simplicial setting.

\begin{Def} 
For an object $Y^\bullet$ in $c\mc{M}$ let {\bf\boldmath $\text{Fib}_sY^\bullet$\unboldmath} denote the fiber of $\Tot_s Y^\bullet\to\Tot_{s-1}Y^\bullet$. 
\end{Def}

\begin{Bem} \label{möglichweise isomorph}
We already mentioned that the spiral exact sequence can be spliced together to an exact couple giving the spectral sequence \Ref{Spiralspektral}:
      $$ E_2^{p,q}=\pi_p\stabhom{X^\bullet}{\Omega^q G}\ \Longrightarrow\ \colim_k\naturalpi{k}{X^\bullet}{\Omega^{p+q-k}G} $$
Like in \cite[3.9]{GoHop:moduli} it follows that this spectral sequence is isomorphic from the $E_2$-term onward to a more familiar one, namely the $G$-cohomology spectral sequence of the total tower $\kl\Tot_sX^\bullet\kr$ for every $G\aus\mc{G}$. Its $E_1$-term consists of
\begin{equation*}
    E_1^s=G^*\text{Fib}_sX^\bullet=\stabhom{\text{Fib}_sX^\bullet}{G}_*.
\end{equation*}
Since $X^\bullet$ is Reedy-fibrant, there is an isomorphism $\text{Fib}_sX^\bullet\cong \Omega^sN^sX^\bullet$ by \cite[p. 391]{GoJar:simp}, where $N^sX^\bullet:=\text{fiber}\hrl X^s\to M^{s}X^\bullet\hrr $ is the geometric normalization of $X^\bullet$. 
Moreover it is true that there is an isomorphism
\begin{equation} \label{G-E1}
    G^*\hrl \text{Fib}_sX^\bullet\hrr =G^*\hrl \Omega^sN^sX^\bullet\hrr \cong N_s\hrl G^{*+s}X^\bullet\hrr ,
\end{equation}
where on the right hand side $N_s$ denotes the normalization of complexes.
Also the spectral sequence differential $d_1:G^*\hrl \text{Fib}_{s+1}X^\bullet\hrr \to G^*\hrl \text{Fib}_sX^\bullet\hrr $ coincides up to sign with the boundary of the normalized cochain complex $N_*\hrl G^*X^\bullet\hrr $. Hence:
\begin{equation} \label{G-cohomology-ss}
    E_2^s = \pi_s\stabhom{X^\bullet}{G}\, \Longrightarrow\, \colim_k\stabhom{\Tot_kX^\bullet}{G} 
\end{equation}
Theorem 6.1(a) from \cite{Boa:ccss} states that the convergence of this spectral sequences is strong if and only if $\Rlim_rE^*_r=0$. 
Problems arise now when we try to relate the target with the term $\Hom{\mc{A}}{\pi^0FX^\bullet}{F_*G}$. Let us formulate all this still in another way. Consider the homology spectral sequence of a cosimplicial space
    $$ E_2^{s,t}=\pi^sF_tX^\bullet\, \Longrightarrow\, \lim_k F_{t-s}\Tot_kX^\bullet $$
with differentials 
    $$ d_r: E_r^{s,t}\to E_r^{s+r,t+r-1} .$$
Again a necessary and sufficient criterion for strong convergence is the vanishing of $\Rlim_rE^*_r$. Now the question arises, what has $\lim_sF\Tot_sX^\bullet$ to do with $F\Tot X^\bullet$.
If our functor $F$ commutes with countable products, we could set up a Milnor-type sequence, which answered this question. But in the cases of interest $F$ will not commute with infinite products, and the question has to remain open in general. 
The point is that $\Tot$ of a fibrant approximation computes a completion of the initial object, which may not coincide with the (localization of the) object. In the case of $F=E_*$ given by a suitable ring spectrum $E$ this will be the $E$-completion, which is proved in \cite{Schorsch:diplom}.
The answer we can offer is that, whenever the convergence results from \cite{Bous:cos-homology-ss} apply or the completion is known to be isomorphic to the localization, we can derive the existence of $X$ in \mc{T}. Later, when it comes to realization questions, we will simply assume that the injective dimension of each object in \mc{A} is finite, which will render all convergence problems trivial.
\end{Bem}

\begin{lemma} \label{Spektralsequenz kollabiert}
Let $X^\bullet$ be an $F$-fibrant object with the property that for $s\uber 0$ 
    $$\pi_s\emph{\stabhom{X^\bullet}{G}} = 0.$$
\begin{punkt}
    \item
Then there is a natural isomorphism
    $$ \lim_sF_*\Tot_s X^\bullet \cong \pi^0F_*X^\bullet .$$
For every $G\aus\{F{\rm -inj}\}$ there are natural isomorphisms
    $$ \colim_k\emph{\naturalpi{k}{X^\bullet}{\Omega^{p+q-k}G}} \cong \colim_s\emph{\stabhom{\Tot_sX^\bullet}{G}}\cong \emph{\Hom{\mc{A}}{\pi^0 FX^\bullet}{FG}}. $$
    \item
If $\pi^0FX^\bullet$ has finite injective dimension then there are isomorphisms
    $$ F\Tot X^\bullet\cong\pi^0FX^\bullet $$
and for every $G\aus\{F{\rm -inj}\}$ 
    $$ \emph{\stabhom{\Tot X^\bullet}{G}}\cong\colim_s\emph{\stabhom{\Tot_s X^\bullet}{G}}\cong\emph{\Hom{\mc{A}}{\pi^0FX^\bullet}{FG}} .$$
\end{punkt}
\end{lemma}

\begin{beweis}
The isomorphisms of part (i) follow from the fact that all the afore\-men\-tioned spectral sequences collapse at the $E_2$-level with non-vanishing terms just in filtration $0$ and converge strongly. Part (ii) follows from ${\rm Fib}_s\Omega^sN^sX^\bullet$ and \Ref{G-E1}, which imply that the constant tower $\{F_*\Tot X^\bullet\}$ and the tower $\{F_*\Tot_kX^\bullet\}$ are pro-equivalent.
\end{beweis}

\begin{kor} \label{Iso in T}
For every $Y$ in \mc{M} let $r^0Y\to Y^\bullet$ be an $F$-fibrant approximation. Then the canonical map $Y\to \Tot Y^\bullet$ is an isomorphism in \mc{T}, if $F$ detects isomorphisms and if $FY$ has finite injective dimension.
\end{kor}

\begin{beweis}
We immediately derive this result from \ref{Spektralsequenz kollabiert}.
\end{beweis}

\begin{Bem} \label{ASS}
Finally we remark that we get back to the modified Adams spectral sequence if we apply the functor $\stabhom{X}{\frei}$ to the total tower of an $F$-fibrant approximation $Y^\bullet$ of an object $Y$ from \mc{T}. The modified Adams spectral sequence is constructed in the same way as the original Adams spectral sequence, but it uses absolute injective resolutions instead of relative ones. It was introduced in \cite{Bri:adams}. Other accounts are given in \cite{Bou:klocal}, \cite{Dev:brown-comenetz} and \cite{Fra:uni}. It can be considered as the Bousfield-Kan spectral sequence of the simplicial space $\map^{\rm pro}\hrl X,Y^\bullet\hrr$. The $E_1$-term is given by
\begin{equation} \label{E_1-exaktes Paar}
     E_1^{s,t} = \pi_0\map\hrl\Sigma^t X, N^sY^\bullet\hrr\cong\stabhom{X}{N^sY^\bullet}_t .
\end{equation}
Since $Y^\bullet$ is an $F$-fibrant approximation to $r^0Y$ it follows that the $E_2$-term takes the following form
    $$ E_2^{s,t} = \Ext{\mc{A}}{s,t}{F_*X}{F_*Y}, $$
which is independent of the choice of the $F$-fibrant approximation and functorial in $X$ and $Y$. The differentials are maps
\begin{equation} \label{ASS-Differential}
     d_r: E_r^{s,t} \to E_r^{s+r,t+r-1}. 
\end{equation}
We have another construction of the modified Adams spectral sequence obtained by applying the functor $\stabhom{\frei}{Y^\bullet}_{F}$ to the Postnikov cotower from \ref{Ko-Postnikov-Turm}. We get an exact couple
\diag{ ...  & \stabhom{\sk_s\hrl X\otimes\Delta^\bullet\hrr}{Y^\bullet}_F \ar[l]\ar[d] & \stabhom{\sk_{s+1}\hrl X\otimes\Delta^\bullet\hrr}{Y^\bullet}_F \ar[l]\ar[d] & ... \ar[l] \\
        & \Ext{\mc{A}}{s+1,s}{A}{B} \ar[ur]_+ & \Ext{\mc{A}}{s+2,s+1}{A}{B} \ar[ur]_+ & }{abgeleitetes Paar}
where the lower terms are identified by \ref{Darstellungseigenschaft von L(N,n)} and \ref{Pushout für ein n-Stück}. The $+$ indicates, that the map raises the external degree by one, so the differentials have the same form as in \Ref{ASS-Differential}. 
It follows from the dual version of \cite[Lemma 3.9.]{GoHop:moduli}, that the exact couple from \Ref{abgeleitetes Paar} is isomorphic to the derived couple of the $E_1$-exact couple in \Ref{E_1-exaktes Paar}.
Hence the spectral sequences coincide from the $E_2$-term on.

For convergence results we have to take the usual precautions, see \cite{Dev:brown-comenetz}. It is shown in \cite{Schorsch:diplom} that it converges strongly to $\stabhom{X}{F^{\wedge}Y}$ if and only if $\Rlim_rE_r^{*,*}=0$. Here $F$ is a topologically flat ring spectrum with $F_*F$ commutative (see \cite{Hov:algebroid}) and $F^{\wedge}Y$ is the $F$-completion of $Y$.
\end{Bem}

\section{Properties of interpolation categories}

In \ref{subsection:Eigenschaften} we describe the properties our tower of interpolation categories enjoys. In subsection \ref{subsection:Moduli} express these facts in terms of moduli spaces. The analogous results their have been obtained in \cite{BlDG:pi-algebra} but our proofs are much shorter due to our truncated model structures.

\subsection{Properties} 
\label{subsection:Eigenschaften}

The axioms that hold for our interpolation categories are taken from \cite[VI.5]{Baues:combinatorial} and briefly explained in appendix \ref{subsection:Ausdehnung}.

\begin{satz} \label{Baues-Turm}
Let $F$ be a homological functor as in \emph{\ref{Annahmen für F}} and $n\bth 1$. 
Then the following diagram
\begin{align*}
    \emph{\Ext{\mc{A}}{n,n}{\pi^0F_*\hrl\frei\hrr}{\pi^0F_*\hrl\frei\hrr}}&\to IP_n\emph{\hrl} F\emph{\hrr} \to IP_{n-1}\emph{\hrl} F\emph{\hrr} \\
           &\to \emph{\Ext{\mc{A}}{n+1,n}{\pi^0F_*\hrl\frei\hrr}{\pi^0F_*\hrl\frei\hrr}} 
\end{align*}
is an exact sequence of categories in the sense of \emph{\ref{exakte Sequenz von Kategorien}}.  
\end{satz}

\begin{beweis}
We have to check the various points in definition \ref{exakte Sequenz von Kategorien}. That the Ext-terms here define natural systems of abelian groups is clear. Property (i) of \ref{exakte Sequenz von Kategorien} is proved in \ref{Ext^(n,n) und Morphismen}.
(ii) follows from \ref{Ext^(n+1,n) und Morphismen}. 
Point (iii) is shown in \ref{Derivation} and (iv) follows from the proof of theorem \ref{Ext^(n+1,n) und Objekte}.  
\end{beweis}

\begin{satz} 
For each $n\bth 0$ the functor $\sigma_n\co IP_{n+1}\emph{\hrl} F\emph{\hrr} \to IP_{n}\emph{\hrl} F\emph{\hrr}$ detects isomorphisms.
\end{satz}

\begin{beweis}
This follows from \ref{Baues-Turm} and \ref{sigma entdeckt Isos}.
\end{beweis}

\begin{satz} \label{Äquivalenz der 0-ten Schicht}
The functor $\pi^0F_*\co IP_0\emph{\hrl} F\emph{\hrr} \to\mc{A}$ is an equivalence of categories.
\end{satz}

\begin{beweis}
We will prove that $\pi^0F$ is essentially surjective and induces a bijection
\begin{equation} \label{Volltreu}
    \Hom{IP_0(F)}{X^\bullet}{Y^\bullet}\to\Hom{\mc{A}}{\pi^0F_*X^\bullet}{\pi^0F_*Y^\bullet} 
\end{equation}
for arbitrary objects $X^\bullet$ and $Y^\bullet$ in $IP_0\hrl F\hrr$. 

Let $A$ be in \mc{A}. Choose an injective resolution $A\to I^\bullet$. 
Using remark \ref{E(I) ist funktoriell} we can realize $I^\bullet$ as a diagram in $\mc{T}=\ho{\mc{M}}$. It was shown in \ref{Existenz von L(A,0)} that we can realize the beginning part of this resolution by a $1$-truncated cosimplicial object $E^\bullet$ in $c_1\mc{M}$. Let $X^\bullet$ in $c\mc{M}$ be $l^1E^\bullet$, where $l^1$ is the left Kan extension to $c\mc{M}$. Such an object is in $IM_0\hrl F\hrr$ because it is an object of the form $L\hrl A,0\hrr$. By construction we have 
    $$\pi^0F_*X^\bullet\cong \ker[\xymatrix@1{ F_*E^0\ar[r]^{d^1-d^0} & F_*E^1}]\cong \ker[I^0\to I^1]\cong A,$$
which proves that $\pi^0F_*$ is essentially surjective. Now let $X^\bullet$ and $Y^\bullet$ be objects in $IP_0\hrl F\hrr$. Suppose we are given a map $A\to B$ in \mc{A} with $A=\pi^0F_*X^\bullet$ and $B=\pi^0F_*B^\bullet$. We can assume that $X^\bullet$ is $0$-$F$-cofibrant and $Y^\bullet$ is $F$-fibrant. Then $X^\bullet$ is of type $L\hrl A,0\hrr$. A map from $A$ to $B$ can be extended to a map
    $$ K\hrl A,0\hrr\to F_*Y^\bullet, $$
since $r^0$ is left adjoint to taking the maximal augmentation $\pi^0\hrl\frei\hrr$. Now \ref{Darstellungseigenschaft von L(N,n)} delivers us a map $L\hrl A,0\hrr=X^\bullet\to Y^\bullet$ in $c\mc{M}$ inducing $A\to B$. Hence the functor is full.

Finally let $X^\bullet\to Y^\bullet$ be a morphism in $IP_0\hrl F\hrr$ that is in the kernel of the map \Ref{Volltreu}. Again we assume, that $X^\bullet$ is $0$-$F$-cofibrant and $Y^\bullet$ is $F$-fibrant. This implies that the morphism is represented by a map $X^\bullet\to Y^\bullet$ in $c\mc{M}$
, but also that $X^\bullet$ is of type $L\hrl A,0\hrr$. 
The induced map
   $$ K\hrl A, 0\hrr \to F_*X^\bullet\to F_*Y^\bullet $$
is null homotopic by assumption. But then the fact that $L\hrl A,0\hrr=X^\bullet\to Y^\bullet$ is null homotopic follows again from \ref{Darstellungseigenschaft von L(N,n)}. This proves that $\pi^0F_*\hrl\frei\hrr$ is faithful.
\end{beweis}

\begin{satz} 
Let $A$ be an object in \mc{A} of injective dimension $\sth n+2$ for $n\bth 0$. The object $A$ is realizable in \mc{T} if and only if there exists an object $X^\bullet$ in $IP_n\emph{\hrl}F\emph{\hrr}$ with $\pi^0F_*X^\bullet\cong A$ or equivalently a potential $n$-stage for $A$. 
\end{satz}

\begin{beweis}
The necessity of the existence of a potential $n$-stage is obvious, but also the sufficiency follows easily, since the obstructions against the existence of a realization as an $\infty$-stage $X^\bullet$ lie in $\Ext{\mc{A}}{n+3+s,n+1+s}{A}{A}$ for $s\bth 0$ and these groups vanish by assumption. Now for a fibrant approximation $X^\bullet\to\wt{X}^\bullet$ the total space $\Tot\wt{X}^\bullet$ is a realization of $A$ by \ref{Spektralsequenz kollabiert}.
\end{beweis}

\begin{satz} 
Let $A$ be an object in \mc{A} with $\dim A\sth n+1$ for $n\bth 0$. Let $X^\bullet$ in $IP_n\emph{\hrl}F\emph{\hrr}$ be an object with $\pi^0F_*X^\bullet\cong A$. Then there exists an object $X$ in \mc{T}, which is a lifting of $X^\bullet$ $($and of $A$$)$, and its isomorphism class in \mc{T} is uniquely determined.
\end{satz}

\begin{beweis}
Analogously to the previous proof now also all obstruction against uniqueness given by theorem \ref{Ext^(n+1,n) und Objekte} vanish.
\end{beweis}

\begin{Def}
Recall that \mc{A} has enough injectives by assumption. 
We consider the full subcategory \boldmath$\mc{T}_n$\unboldmath\ of \mc{T} consisting of those objects $X$ such that $F_*X$ has injective dimension $\sth n$. This defines an increasing filtration of \mc{T} with $\mc{T}_0$ equal to the full subcategory of $F$-injective objects $\mc{T}_{\rm inj}$. The inclusion functors $\mc{T}_n\inj\mc{T}$ will be called \boldmath$i_n$\unboldmath.
\end{Def}

\begin{satz} \label{T_{n}=True_{n}(F)}
The functors $\theta_{k} i_n:\mc{T}_n\to IP_{k}\emph{\hrl}F\emph{\hrr}$ are full for $k\bth n-1$. The functors $\theta_{k} i_n:\mc{T}_n\to IP_{k}\emph{\hrl}F\emph{\hrr}$ are faithful for $k\bth n$. 
\end{satz}


\begin{beweis}
If $X$ is an object of $\mc{T}$ then its image $\theta_kX$ in $IP_k\hrl F\hrr$ is the $k$-$F$-equivalence class of $r^0X$. Let $X$ and $Y$ be in $\mc{T}_{n}$ where we assume from the beginning on that both are fibrant and cofibrant. We have to show that the map
\begin{equation} 
     \Hom{\mc{T}_{n}}{X}{Y}\to\Hom{IP_{k}(F)}{\theta_{k}X}{\theta_{k}Y}
\end{equation}
is a bijection for $k\bth n-1$. To prove surjectivity we take $F$-fibrant replacements $\wt{X}^\bullet$ and $Y^\bullet$ of $\theta_{k}X=r^0X$ and $\theta_{k}Y=r^0Y$ respectively and then we replace $\wt{X}^\bullet$ Reedy cofibrantly by $X^\bullet$. Now each morphism $\hel f\her$ in $\Hom{IP_{k}(F)}{r^0X}{r^0Y}\cong \Hom{IP_{k}(F)}{X^\bullet}{Y^\bullet}$ is represented by a map 
    $$ f:\sk_{k+1}X^\bullet\to Y^\bullet $$
in $c\mc{M}$. The obstructions against extending this map to higher skeleta of $X^\bullet$ lie in $\Ext{\mc{A}}{k+2+s,k+1+s}{F_*X}{F_*Y}$ for $s\bth 0$ by \ref{Ext^(n+1,n) und Morphismen}. All these groups vanish for $k\bth n-1$ because the injective dimension is smaller than or equal to $n\unter k+2$. We end up with a map $f_\infty:X^\bullet\to Y^\bullet$. Now we get a morphism 
    $$ \xymatrix{\wt{f}:X \cong \Tot X^\bullet\ar[r]^-{\Tot f_\infty} & \Tot Y^\bullet\cong Y }$$
in \mc{T}, where the isomorphisms are the canonical maps from \ref{Iso in T}, and they are isomorphisms since $F$ detects them. By lemma \ref{Tot-Delta-Quillen-Paar} or remark \ref{Dasselbe Quillenpaar} $R\Tot$ and $Lr^0$ are a Quillen pair, and so $\sigma_n\wt{f}$ corresponds to $[f]$ via the isomorphism
   $$\Hom{IP_{k}(F)}{r^0X}{r^0Y}\cong \Hom{IP_{k}(F)}{X^\bullet}{Y^\bullet}$$
induced by the various replacements. So we have shown that $\theta_n$ is full.

The second part of the theorem amounts to prove the injectivity of the map
\begin{equation} 
     \Hom{\mc{T}_{k}}{X}{Y}\to\Hom{IP_{k}(F)}{\theta_{k}X}{\theta_{k}Y}
\end{equation}
for $k\bth n$. This map is a homomorphism of abelian groups since $\theta_k$ is additive. Let $g:X\to Y$ represent a morphism that is mapped to zero. 
Again we pick replacements $X^\bullet$ and $Y^\bullet$ of $r^0X$ and $r^0Y$ as above. We find a map $g_\infty:X^\bullet\to Y^\bullet$ whose homotopy class is uniquely determined by $r^0g:r^0X\to r^0Y$ and which is nullhomotopic in $c\mc{M}^F$ when we restrict it to the $\hrl k+1\hrr$-skeleton of $X^\bullet$.
This is displayed in the following solid arrow diagram which strictly commutes:
\diagr{ \sk_{k+1}X^\bullet \ar[r]^-{H}\ar[d]_{s_{k+1}}\ar[dr]^->>>{j_{k+1}} & \hom\hrl\Delta^1,Y^\bullet\hrr \ar[r]^-{d_0} \ar[dr]_-<{d_1} & Y^\bullet \\
        \sk_{k+2}X^\bullet \ar[r]_-{j_{k+2}}\ar@{.>}[ur]|{\phantom{xx}}^-<<{H'} & X^\bullet \ar[r]_\ast\ar[ur]|{\phantom{xx}}_->>>{g_\infty} & Y^\bullet}
The evaluation maps $d_0$ and $d_1:\hom\hrl\Delta^1,Y^\bullet\hrr\to Y^\bullet$ are $F$-equivalences, so for both objects their $\pi^0F_*$-term is isomorphic to $F_*Y$ in $\mc{T}_n$. In particular it follows from the first part of the theorem that $H'$ exists with $H's_{k+1}\simeq H$ in the $F$-structure. 
Actually the proof of \ref{Ext^(n+1,n) und Morphismen} shows that we can arrange this to be strictly equal. It tells us that $d_0H'$ and $g_\infty j_{k+2}$ are both extensions of the map $g_{\infty}j_{k+1}=d_0H$. The obstructions against uniqueness of liftings, which are the homotopy classes of these extensions, lie in $\Ext{\mc{A}}{k+1,k+1}{F_*X}{F_*Y}$ and this group vanishes since the injective dimension is smaller than or equal to $n\unter k+1$ by assumption. It follows that $g_{\infty}j_{k+2}$ is $F$-homotopic to $d_0H'$. The same argument works with the other evaluation map $d_1$ and shows that $g_{\infty}j_{k+2}$ is nullhomotopic. By induction we can extend this over all skeleta since all higher obstruction groups also vanish. 
The skeletal tower of $X^\bullet$ is a tower of $F$-cofibrations between Reedy cofibrant (aka. $F$-cofibrant) objects since $X^\bullet$ is Reedy cofibrant, therefore we have
    $$ X^\bullet\cong\colim_k\sk_{k}X^\bullet\cong\hocolim_k\sk_kX^\bullet. $$
Hence the successive extensions give us a map $g_{\infty}':X^\bullet\to Y^\bullet$ which on one side is homotopic to $g_{\infty}$ and on the other to $*$. 
Because the homotopy class of $g_{\infty}$ or $g_{\infty}'$ corresponds under the isomorphism
     $$ \pi_0\map\hrl X^\bullet,Y^\bullet\hrr\cong\pi_0\map\hrl r^0X,r^0Y\hrr$$
to $r^0g$, this shows that our original map $r^0g$ is nullhomotopic in $c\mc{M}$. Finally constant cosimplicial objects over fibrant objects are Reedy fibrant, so we can apply and conclude
    $$ [g]=0\ \aus\ \pi_0\map_{\mc{M}}\hrl X,Y\hrr, $$
since $\Tot$ maps external homotopies between Reedy fibrant objects to homotopies in \mc{M}. 
\end{beweis}

Actually the previous statement can be strengthened since for both assertions only the fact that $Y$ is in $\mc{T}_n$ was needed. A plausible extension of \ref{T_{n}=True_{n}(F)} is the statement, that the functor $\mc{T}\to IP_n\hrl F\hrr$ given by
    $$ X \mapsto \Hom{IP_n(F)}{\theta_n X}{Y^\bullet} $$
for some $Y^\bullet$ in $IP_n\hrl F\hrr$ is representable if and only if $\dim\pi^0FY^\bullet\sth n$, but we have not been able to prove that.

The following theorem relates the tower of interpolation categories to the Adams spectral sequence associated to $F$.
\begin{Def}
A map $f$ in \mc{T} is said to be of {\bf\boldmath Adams filtration $n$} if it admits a factorization $f=f_1...f_n$, where the maps $f_i$ induce the zero map via $F$ in \mc{A}. Let $F^n\stabhom{X}{Y}$ be the set of all maps of Adams filtration $n$, where we set $F^0\stabhom{X}{Y}:=\stabhom{X}{Y}$.
We obtain a decreasing filtration of $\stabhom{X}{Y}$.
\end{Def}

\begin{satz} \label{Adams-Filtrierung}
For $n\bth 0$ there is a natural isomorphism
    $$ F^{n+1}\emph{\stabhom{X}{Y}} \cong\ker\emph{\hel\stabhom{X}{Y}}\to\emph{\Hom{IP_n(F)}{\theta_nX}{\theta_nY}\her}.$$
\end{satz}

\begin{beweis}
Let $r^0Y\to Y^\bullet$ be an $F$-fibrant approximation and remember that $NF_*Y^\bullet$ is an injective resolution of $F_*Y$.
For $n=0$ the claim follows from the equivalence $IP_0\hrl F\hrr\cong\mc{A}$. For $n=1$ a map $f\aus F^2$ in particular induces the zero map in \mc{A}, thus it admits a factorization
     $$ X\to \wt{Y}^1\to Y $$
where $\wt{Y}^1$ is the fiber of the map $Y\to Y^0$. It follows easily that such maps $f$ are characterized by the fact that the map $X\to\wt{Y}^1$ induces the zero map in \mc{A}. On the other hand a map $f\co X\to Y$ is in the kernel above if and only if the map $\sk_{2}\hrl X\otimes\Delta^\bullet\hrr\to Y^\bullet$ representing $\theta_1f\co\theta_1X\to\theta_1Y$ admits a lifting:
\diagr{  & \hom\hrl\Delta^1,Y^\bullet\hrr\times_{Y^\bullet}\ast \ar[d] \\
        \sk_{2}\hrl X\otimes\Delta^\bullet\hrr \ar[ru]\ar[r] & Y^\bullet }
Considering the non-degenerate $1$-simplex in $X\otimes\Delta^1$ shows that the map $X\to Y^1$ is null homotopic. This map also factorizes over $X\to\wt{Y}^1$ and, since $F\wt{Y}^1\to F_*Y^1$ is injective, it follows that $X\to\wt{Y}^1$ induces the zero map.
We can work backwards and show that, if $X\to\wt{Y}^1$ induces the zero map, we can construct null homotopy on the $2$-skeleton, which proves the isomorphism.
For higher $n$ we proceed inductively. We show that the map $X\to\wt{Y}^n$, where $\wt{Y}^n$ is the fiber of the map $\wt{Y}^{n-1}\to Y^{n-1}$, induces the zero map if and only if there is a diagram:
\diagr{  & \hom\hrl\Delta^1,Y^\bullet\hrr\times_{Y^\bullet}\ast \ar[d] \\
        \sk_{n+1}\hrl X\otimes\Delta^\bullet\hrr \ar[ru]\ar[r] & Y^\bullet }
\end{beweis}

\subsection{Moduli spaces of realizations}
\label{subsection:Moduli}

Again we let $F$ be a homological functor as in \ref{Annahmen für F}, $F$ {\bf detects isomorphisms} in \mc{T}. 
Note that this assumption means that a map $X\to Y$ in \mc{M} induces an isomorphism $F_*X\to F_*Y$ if and only if it was a weak equivalence. We summarize the necessary theory of moduli spaces in \ref{appendix:moduli}.

\begin{Def}
Let $A$ be an object in \mc{A}. We define the {\bf space of realizations} of $A$ to be the moduli space of all objects $X$ in \mc{M}, such that their image $F_*X$ is isomorphic to $A$ (see Def. \ref{Modulraum eines Objektes}). We will write \boldmath${\bf\rm Real}\hrl A\hrr$\unboldmath.

We define the {\bf space of \boldmath$n$\unboldmath-th partial realizations} of $A$ to be the moduli space of all objects $X^\bullet$ in $c\mc{M}$ that are potential $n$-stages for $A$ (see Def. \ref{potentielles n-Stück}). We will write \boldmath${\bf\rm Real}_n\hrl A\hrr$\unboldmath. Everything makes also sense for $n=\infty$ and hence we define in the same way the {\bf space of \boldmath$\infty$\unboldmath-stages} of $A$ and denote it by \boldmath${\bf\rm Real}_\infty\hrl A\hrr$\unboldmath.
Recall that $\infty$-\mc{G}-structure is just another name for the \mc{G}-structure. 
\end{Def}

The first theorem we are heading for is \ref{Unendlich-Stücke sind Realisierungen} which tells us that $\infty$-stages are the same as actual realizations in $\mc{T}=\ho{\mc{M}}$. The next step is theorem \ref{holim der Modulräume} which relates the moduli space of $\infty$-stages to the spaces $\Real_n\hrl A\hrr$ of potential $n$-stages.
Finally we establish in \ref{Modulraumsequenz für Objekte} a fiber sequence involving $\Real_{n-1}\hrl A\hrr$ and $\Real_n\hrl A\hrr$.

\begin{Bem} \label{die Abbildung Real_8(A) -> Real(A)}
To relate an $\infty$-stage of an object in \mc{A} to an actual realization we use the functor $\Tot:c\mc{M}\to\mc{M}$. By \ref{möglichweise isomorph} there is a spectral sequence:
    $$ E_2^{s,t} = \pi_s F_tX^\bullet\, \Longrightarrow\, \lim_k F_{t-s}\Tot_kX^\bullet$$
From lemma \ref{Spektralsequenz kollabiert} we can read off that for an $\infty$-stage $X^\bullet$ of an object $A$ with finite injective dimension the spectral sequence collapses and its edge homomorphism gives an isomorphism
    $$ F_*\Tot X^\bullet\cong\lim_kF_*\Tot_kX^\bullet\cong A. $$
More generally the spectral sequence gives such an isomorphism whenever the results in \cite{Bous:cos-homology-ss} say so.
We see that under these assumptions the functor $\Tot$ induces a natural map
\begin{equation} \label{von Unendlichstück zu Realisierung}
    \Real_{\infty}\hrl A\hrr \to \Real\hrl A\hrr .
\end{equation}
\end{Bem}

\begin{satz} \label{Unendlich-Stücke sind Realisierungen}
The map \emph{\Ref{von Unendlichstück zu Realisierung}} is a weak equivalence of spaces if $A$ has finite injective dimension or if the convergence results from \emph{\cite{Bous:cos-homology-ss}} apply.
\end{satz}

\begin{beweis}
Let $X$ be a realization of $A$ in \mc{M}. 
Then the canonical map $X\to \Tot r^0X=\Tot_0r^0X=X$ is even an isomorphism in \mc{M}. 

Let $X^\bullet$ be a vertex in $\Real_{\infty}\hrl A\hrr$, in other words an $\infty$-stage of $A$. Without loss of generality we assume that $X^\bullet$ is $F$-fibrant and Reedy cofibrant, because these manipulations induce self equivalences of the moduli space $\Real_{\infty}\hrl A\hrr$. But now the map 
    $$ r^0\Tot X^\bullet \to X^\bullet $$
is an $F$-equivalence by \ref{Spektralsequenz kollabiert}. Since $F$ detects isomorphisms in \mc{T}, this shows that the maps induced by $\Tot$ and $r^0$ are mutually inverse homotopy equivalences.
\end{beweis}

The rest of this subsection is true without any restriction on $A$.
\begin{satz} \label{holim der Modulräume}
The canonical map 
    $$ \Real_{\infty}\emph{\hrl} A\emph{\hrr}\to \holim_n \Real_n\emph{\hrl} A\emph{\hrr} $$
is a weak equivalence. 
\end{satz}

We prove this theorem after having established two lemmas.
\begin{Def}
Let \boldmath${\rm\bf weak}_S\hrl A^\bullet,B^\bullet\hrr$\unboldmath\ denote the simplicial set given by
    $${\rm weak}_S\hrl A^\bullet,B^\bullet\hrr_n:= \Hom{\mc{W}_S}{A^\bullet\otimes\Delta^n}{B^\bullet} ,$$
where $\mc{W}_S$ is the subcategory of weak equivalences in some simplicial model structure $S$ on $c\mc{M}$. If $A^\bullet$ is fibrant and cofibrant in $S$ then 
    $${\rm weak}_S\hrl A^\bullet,A^\bullet\hrr = \text{haut}_S\hrl A^\bullet\hrr$$
by definition \ref{haut}. 
Analogously to remark \ref{haut als Komponenten von map} we observe that ${\rm weak}_S\hrl A^\bullet, B^\bullet\hrr$ is a union of connected components of $\map_S\hrl A^\bullet, B^\bullet\hrr$.
\end{Def}

\begin{lemma} \label{weak und Limites}
Let \mc{G} be a class of injective models for \mc{M}. Let $X^\bullet$ be a Reedy cofibrant object and $Y^\bullet$ be a \mc{G}-fibrant object in $c\mc{M}$. Then there is a canonical map
\begin{align*}
    \holim_n{\rm weak}_{n-\mc{G}}\emph{\hrl}\sk_{n+1}X^\bullet, Y^\bullet\emph{\hrr} &\stackrel{\simeq}{\longrightarrow} \lim_n{\rm weak}_{n-\mc{G}}\emph{\hrl}\sk_{n+1}X^\bullet, Y^\bullet\emph{\hrr} \\
        & \cong {\rm weak}_{\mc{G}}\emph{\hrl} X^\bullet, Y^\bullet\emph{\hrr} 
\end{align*}
where the first one is a weak equivalence and the second one is an isomorphism which are natural in both variables for $\mc{G}$-equivalences.
\end{lemma}

\begin{beweis}
First we observe that the corresponding statement for the functor $\map\hrl\frei,\frei\hrr$ is true. Here $\map\hrl\frei,\frei\hrr$ which is the external mapping space from \ref{Externe simpliziale Struktur} always has homotopy meaning since $\sk_{n+1}X^\bullet\to X^\bullet$ is an $n$-\mc{G}-cofibrant approximation. Also the tower maps are fibrations by (SM7') because they are induced by the \mc{G}-cofibration $\sk_nX^\bullet\to\sk_{n+1}X^\bullet$. Finally $\map\hrl\frei,\frei\hrr$ turns colimits in the first variable into limits and $\colim_n\sk_{n+1}X^\bullet\cong X^\bullet$. 

The proof is finished by the above remark that ${\rm weak}_S\hrl X^\bullet, Y^\bullet\hrr$ is a union of components in $\map_S\hrl X^\bullet, Y^\bullet\hrr$ and that these components form a tower because $n$-\mc{G}-equivalences are mapped to $(n\!-\!1)$-equivalences by the restriction of the upper maps.
\end{beweis}

\begin{lemma} \label{haut ist Limes}
Let $X^\bullet$ be $F$-fibrant and Reedy cofibrant, then the canonical map
    $$ {\rm haut}_F\emph{\hrl}X^\bullet\emph{\hrr} \to \holim_n {\rm haut}_{n-F}\emph{\hrl}\sk_{n+1} X^\bullet\emph{\hrr} $$
is a weak equivalence. 
\end{lemma}

\begin{Bem}
Note that the homotopy self equivalences on the right hand side can also be taken in the $F$-structure since $n$-$F$-equivalences and $F$-equivalences agree on $n$-$F$-cofibrant objects.
\end{Bem}

\begin{proofof}{\ref{haut ist Limes}}
The inclusions of the skeletons into $X^\bullet$ induce the following commutative diagram:
\diagr{\text{haut}_{n-F}\hrl \sk_{n+1}X^\bullet\hrr \ar[d]_{\simeq}\ar[r] & \text{haut}_{(n-1)-F}\hrl \sk_{n}X^\bullet\hrr \ar[d]^{\simeq} \\
      {\rm weak}_{n-F}\hrl\sk_{n+1}X^\bullet,X^\bullet\hrr \ar[r]^-{\omega_n} & {\rm weak}_{(n-1)-F}\hrl\sk_{n}X^\bullet,X^\bullet\hrr }
Both horizontal maps fit into a tower of maps when we vary $n$. We want to compute the homotopy limit of the upper tower. To do this we have to replace this tower by an objectwise weakly equivalent one in which the tower maps are fibrations. This is provided by the lower tower as we proved in \ref{weak und Limites}. 
The vertical maps are homotopy equivalences because $\sk_{n+1}X^\bullet\to X^\bullet$ is a cofibrant approximation in the $n$-$F$-structure. 
The result follows now from \ref{weak und Limites}.
\end{proofof}

\begin{proofof}{\ref{holim der Modulräume}}
By theorem \ref{DK-Charakterisierung von Modulraeumen} we have the following weak equivalences
\begin{equation*}
    \Real_{\infty}\hrl A\hrr  \simeq \bigsqcup_{\langle X^\bullet\rangle_F} B\text{haut}_F\hrl X^\bullet\hrr 
\end{equation*}
where the coproduct is taken over all $F$-equivalence classes $\langle X^\bullet\rangle$ of $\infty$-stages $X^\bullet$ of $A$. By the same theorem we obtain the first of the next two weak equivalences
\begin{equation*}
    \Real_n\hrl A\hrr \simeq \bigsqcup_{\langle X^\bullet_n\rangle_F} B\text{haut}_{F}\hrl X^\bullet_n\hrr \simeq \bigsqcup_{\langle X^\bullet_n\rangle_{n-F}} B\text{haut}_{n-F}\hrl X^\bullet_n\hrr ,
\end{equation*}
where the coproduct is taken over all $F$-equivalence classes $\langle X^\bullet_n\rangle$ of potential $n$-stages $X^\bullet_n$ of $A$. Because the $F$- and the $n$-$F$-equivalences agree on $n$-$F$-cofibrant objects, there is a one-to-one correspondence between equivalence classes of potential $n$-stages in the $F$-structure and in the $n$-$F$-structure and there is a weak equivalence ${\rm haut}_F\hrl X^\bullet\hrr\simeq{\rm haut}_{n-F}\hrl X^\bullet\hrr$. Hence we get the second weak equivalence where the coproduct is taken over weak equivalence classes in the $n$-$F$-structure.
The theorem follows now from lemma \ref{haut ist Limes} and the fact that the classifying space functor $B$ from simplicial monoids to \mc{S} preserves weak equivalences, fibrations and limits.
\end{proofof}

\begin{satz} \label{0-ter Modulraum für Objekte}
Let $A$ be an object in \mc{A}. Then we have a weak equivalence
    $$ \Real_0\emph{\hrl} A\emph{\hrr}\simeq B{\rm Aut}\emph{\hrl} A\emph{\hrr} .$$
\end{satz}

\begin{beweis}
This is just a restatement of the equivalence in \ref{Modulraum von L(A,0)}.
\end{beweis}

\begin{satz} \label{Modulraumsequenz für Objekte}
Let $X^\bullet_{n-1}$ be a potential $\emph{\hrl} n-1\emph{\hrr}$-stage for an object $A$ in \mc{A}. Then there is a fiber sequence
    $$ \mc{H}^{n+1}\emph{\hrl} r^0A,A\emph{\hel} n\emph{\her}\emph{\hrr} \to \Real_n\emph{\hrl} A\emph{\hrr}_{X^\bullet_{n-1}} \to \mc{M}\emph{\hrl} X^\bullet_{n-1}\emph{\hrr},  $$
where $\Real_n\emph{\hrl}A\emph{\hrr}_{X^\bullet_{n-1}}$ are those components of $\Real_n\emph{\hrl} A\emph{\hrr}$ that correspond to objects $X^\bullet$ with $\sk_nX^\bullet\simeq X^\bullet_{n-1}$.
\end{satz}

\begin{beweis}
By \ref{Pushout für ein n-Stück} there is a cofiber sequence
    $$ X^\bullet_{n-1}\to X^\bullet_n\to L\hrl A\hel n\her,n\hrr $$
in $c\mc{M}^F$ inducing the following fiber sequence in \mc{S}:
    $$ \map\hrl L\hrl A\hel n\her,n\hrr, X^\bullet_n\hrr \to \map\hrl X^\bullet_n, X^\bullet_n\hrr \to \map\hrl X^\bullet_{n-1}, X^\bullet_n\hrr $$
Passing to appropriate components gives a fiber sequence
    $$ \map\hrl L\hrl A\hel n\her,n\hrr, X^\bullet_n\hrr \to {\rm weak}_{n-F}\hrl X^\bullet_n, X^\bullet_n\hrr \to {\rm weak}_{(n-1)-F}\hrl X^\bullet_{n-1}, X^\bullet_n\hrr $$
of grouplike simplicial monoids. Applying the classifying space functor $B$ to this sequence yields a fiber sequence
    $$ B\ \map\hrl L\hrl A\hel n\her,n\hrr, X^\bullet_n\hrr \to \mc{M}\hrl X^\bullet_n\hrr_{X^\bullet_{n-1}} \to \mc{M}\hrl X^\bullet_{n-1}\hrr. $$
Let $\Gamma: \text{CoCh}^{\bth 0}\hrl\mc{A}\hrr\to c\mc{A}$ be the Dold-Kan functor. We compute finally using \ref{Darstellungseigenschaft von L(N,n)}:
\begin{align*}
    B\ \map\hrl L\hrl A\hel n\her,n\hrr, X^\bullet_n\hrr & \simeq B\ \map\hrl K\hrl A\hel n\her,n\hrr, \sk_{n+1}F_*X^\bullet_n\hrr \\
        & \simeq B\ \map\hrl K\hrl A\hel n\her,n\hrr, r^0A\hrr \\
        & \simeq B\ \Gamma \hrl\Hom{\mc{A}}{A}{A}\hel n\her_{\rm ext}\hrr \\
        & \simeq \Gamma \hrl\Hom{\mc{A}}{A}{A}\hel n+1\her_{\rm ext}\hrr \\
        & \simeq \mc{H}^{n+1}\hrl A,A\hel n\her\hrr,
\end{align*}
where $\Hom{\mc{A}}{A}{A}\hel n\her_{\rm ext}$ is viewed as a cochain complex concentrated in degree $n$. Here $\hel 1\her_{\rm ext}$ is the external shift from \ref{mehr Eigenschaften von W und W-quer}.
\end{beweis}

\begin{satz} \label{0-ter Modulraum für Morphismen}
Let $f$ be a morphism in \mc{A}. Then we have a weak equivalence
    $$ \Real_0\emph{\hrl} f\emph{\hrr}\simeq B{\rm Aut}\emph{\hrl} f\emph{\hrr} .$$
\end{satz}

\begin{beweis}
This follows readily from the equivalence $\mc{A}\cong IP_0\hrl F\hrr$ of categories of \ref{Äquivalenz der 0-ten Schicht}.
\end{beweis}

\begin{satz} \label{Modulraumsequenz für Morphismen}
Let $f:X^\bullet_n\to Y^\bullet_n$ be a map of potential $n$-stages for objects $A$ and $B$ respectively in \mc{A}. 
Then there is a fiber sequence
    $$ \mc{H}^{n}\emph{\hrl} A\emph{\hel} n\emph{\her},B\emph{\hrr} \to \mc{M}\emph{\hrl} f\emph{\hrr}_{\sk_n f} \to \mc{M}\emph{\hrl} \sk_n f\emph{\hrr}  $$
\end{satz}

\begin{beweis}
We can assume without loss of generality that $X^\bullet_n$ and $Y^\bullet_n$ are Reedy cofibrant and $F$-fibrant.
As in the proof of \ref{Modulraumsequenz für Objekte} we obtain a fiber sequence
    $$ \map\hrl L\hrl A\hel n\her,n\hrr, Y^\bullet_n\hrr \to \map\hrl X^\bullet_n, Y^\bullet_n\hrr \to \map\hrl\sk_nX^\bullet_n, Y^\bullet_n\hrr .$$ 
Proceeding like in the previous proof we arrive at the conclusion.
\end{beweis}

\section{Examples and applications}
\label{section:Anwendungen}

All the applications given here will just involve the obstruction calculus. We do not have yet applications of the interpolation categories themselves.
\subsection{Very low dimensions}

\begin{Beisp}
If the injective dimension of the target category is $0$ then the tower of interpolation categories simply collapses to the equivalences:
\diagr{ IP_0\hrl F\hrr \ar[r]^-{\cong} & \mc{A} \\ 
       \mc{T}=\mc{T}_0 \ar[u]^{\cong}\ar[r]_-{\cong}^-F & \mc{A}_{\rm inj}\ar[u]_{\cong}  }
Here the lower equivalence was already stated in \ref{E(I) ist funktoriell}. 
\end{Beisp}

\begin{Beisp}
If the injective dimension of \mc{A} is $1$ then the tower of interpolation categories has one non-trivial step:
\diagr{ \mc{T}=\mc{T}_1 \ar[r]^{\cong} & IP_1\hrl F\hrr \ar[d] &  \\ 
        & IP_0\hrl F\hrr \ar[r]^-{\cong} & \mc{A} \\ 
        & \mc{T}_0 \ar[u]\ar[r]_-{\cong} & \mc{A}_{\rm inj}\ar[u]  }
We can express this using \ref{Ext^(n,n) und Morphismen} or \ref{Baues-Turm} by saying that
    $$ \Ext{\mc{A}}{1,1}{F_*\hrl\frei\hrr}{F_*\hrl\frei\hrr} \to \mc{T} \to \mc{A} $$
is a linear extension of categories which is defined in \cite[VI.5]{Baues:combinatorial}.
\end{Beisp}

\newcommand{\Kp}{\ensuremath{{K_{(p)}}\mbox{}}}%
\newcommand{\En}{\ensuremath{{E(n)}\mbox{}}}%
\newcommand{\Pic}{\ensuremath{{\rm Pic}\mbox{}}}%
\subsection{Some $\En$-local Picard groups}
\label{subsection:Picard}

Let \mc{C} be a symmetric monoidal category whose pairing is called smash product and denoted by $\wedge$. An invertible object $X$ is one such that there exists a $Y$ in \mc{C} with $\sm{X}{Y}\cong S$, where $S$ denotes the unit of the monoidal structure. The isomorphism classes of invertible objects inherit an abelian group structure which we will call the Picard group $\Pic(\mc{C})$. It is an abelian group, but sometimes in a higher universe. It was defined by Hopkins and we refer to \cite{Hop-Mah-Sad:Picard}, where it is proved, that the Picard group of the whole stable homotopy category of spectra is $\mathbbm{Z}$. There are also computations involving the Picard group of the $K(n)$-local category, where $K(n)$ denotes $n$-th Morava $K$-theory. 

Fix a prime $p$.
The problem of computing $\Pic(\En)$ for the $\En$-local category of spectra was considered in \cite{Hov-Sad:invertible}. Here $\En$ denotes the $n$-th Johnson-Wilson spectrum. It is a Landweber exact theory with 
    $$\En_*=\mathbbm{Z}_{(p)}[v_1,...,v_{n-1},v_n^{\pm 1}]$$
where $|v_i|=2(p^i-1)$. $E(1)$ is a retract of $\Kp$, which is sometimes called the Adams summand. Hovey and Sadofsky prove that there is a splitting
   $$ \Pic\hrl{E\hrl n\hrr_*}\hrr\cong\mathbbm{Z}\oplus\Pic^0\hrl{E\hrl n\hrr_*}\hrr, $$
where $\mathbbm{Z}$ is generated by the unit of the smash product, the stable localized sphere $L_nS^0$, and they show
   $$ \Pic^0\hrl{E\hrl n\hrr}\hrr=0 ,$$
whenever $2p-1\uber n^2+n$. Their argument is the following: First they prove in \cite[2.4.]{Hov-Sad:invertible} that for an element $X\aus\Pic^0\hrl\En\hrr$ there is an isomorphism
\begin{equation}\label{En-iso}
    \En_\ast X\cong\En_{\ast}
\end{equation}
as $\En_\ast \En$-comodules. This turns the question into a moduli problem adres\-sed in this paper. Actually the Picard group in the presence of the isomorphism \Ref{En-iso} is nothing but $\pi_0\Real(\En_*)$.
Then they show in \cite[5.1.]{Hov-Sad:invertible}, that for $p\uber n+1$ the category of $\En_\ast \En$-comodules has injective dimension $\sth n^2+n$. This is also proved in \cite[Theorem 9]{Fra:uni}.
Considering the $\En$-Adams spectral sequence, which they prove to converge nicely, they see, that the first obstruction for realizing the isomorphism \Ref{En-iso} as a map $X\to L_nS^0$ lies in $\Ext{\En_\ast \En}{2p-1,2p-2}{E\hrl n\hrr_*}{E\hrl n\hrr_*}$ by the usual sparseness in the chromatic setting. Now the statement is clear, since for $2p-1\uber n^2+n$ this obstruction group is zero. 
Using their vanishing line we can extend the range of calculations of Picard groups slightly.

\begin{satz} \label{Pic=Ext}
Fix a prime $p$ and a natural number $n$ such that $p\uber n+1$ and $4p-3\uber n^2+n$. Then we have:
    $$ \Pic^0\emph{\hrl} E\emph{\hrl} n\emph{\hrr}\emph{\hrr}\cong\emph{\Ext{E(n)_*E(n)}{2p-1,2p-2}{E\hrl n\hrr_*}{E\hrl n\hrr_*}} .$$
\end{satz}

\begin{beweis}
According to \Ref{En-iso} we are trying to realize the object $\En_*$. In the range under consideration the only obstruction against existence lies in 
    $$\Ext{\En_*\En}{2p,2p-2}{E\hrl n\hrr_*}{E\hrl n\hrr_*}, $$ 
but the obstruction vanishes because we know that there is a realization, $L_nS^0$. Now theorem \ref{Ext^(n+1,n) und Objekte} tells us that  the uniqueness obstruction group 
    $$\Ext{\En_*\En}{2p-1,2p-2}{E\hrl n\hrr_*}{E\hrl n\hrr_*}$$ 
acts freely and transitively on the realizations. All these elements give actual spectra, since all other obstruction groups in the asserted range vanish by \cite[5.1.]{Hov-Sad:invertible}.
\end{beweis}

For the range $2p-1\uber n^2+n$ we get back the result of Hovey and Sadofsky. \\

{\bf\large Acknowledgements}\\

This article is the second part of my thesis carries out at the Universit\"at Bonn.
I would like to thank my advisor Jens Franke and my coadvisor Stefan Schwede for their interest and help over the years. I also would like to thank Paul Goerss for his encouragement.
Finally I am deeply indebted to Pete Bousfield for providing me with a copy of his then unpublished preprint \cite{Bou:cos} and for several extremely enlightening e-mails.
I was supported by the Bonner Internationale Graduiertenschule (BIGS) which gave me the opportunity to visit several conferences on homotopy theory.

\appendix
\section{The external simplicial structure on $c\mc{M}$}
\label{appendix:external}
The resolution model structures are not compatible with the internal simplicial structure. Here we describe the external simplicial structure, which will be compatible with the resolution structure and its truncated versions.

\begin{Bem} \label{Tensor als Koende}
For $X^\bullet$ in $c\mc{M}$ and $L$ in \mc{S} we can perform the following coend construction: Let $\bigsqcup_{L_\ell}X^m$ be the coproduct in \mc{M} of copies of $X^m$ indexed by the set $L_\ell$, and view this as a functor $\Delta^{\rm op}\times\Delta\to\mc{M}$. Then we can take the coend
    $$ X^\bullet\otimes_{\Delta}L:=\int^{\Delta}\bigsqcup_{L_\ell}X^m \aus\mc{M} .$$
\end{Bem}

We are now ready to describe the functors that will enrich all our model structures to simplicial model categories.

\begin{Def} \label{Externe simpliziale Struktur}
We define a simplicial structure on $c\mc{M}$. Let $K$ be in \mc{S} and $X^\bullet$ and $Y^\bullet$ in $c\mc{M}$, then set
\begin{equation*}
          \hrl X^\bullet\otimes^{\rm ext} K\hrr^n  :=\, X^\bullet\otimes_{\Delta} \hrl K\times\Delta^n\hrr , 
\end{equation*}
where $\times$ denotes the usual product of simplicial sets and $\Delta^n$ is the standard $n$-simplex,
\begin{equation*}
      \hspace*{-2.2cm} \hom^{\rm ext}\hrl K,X^\bullet\hrr^n :=\ \ \prod_{K_n}X^n, 
\end{equation*} 
where the product is taken over the set of $n$-simplices of $K$, and finally
\begin{equation*}
      \hspace*{5mm}    \map^{\rm ext}\hrl X^\bullet,Y^\bullet\hrr_n :=\ \Hom{c\mc{M}}{
X^\bullet\otimes^{\rm ext}\Delta^n}{Y^\bullet}.
\end{equation*} 
We call this the {\bf external (simplicial) structure} on $c\mc{M}$. Note that we do not refer to any simplicial structure of \mc{M}. 
We will usually drop the superscripts.
\end{Def} 

\begin{Def} \label{äußere Einhängung}\label{aeussere Einhaengung}
For an object $X^\bullet$ in $c\mc{M}$ we define its {\bf\boldmath $s$-th external suspension\unboldmath} $\Sigma_{\rm ext}^sX^\bullet$ by the following pushout diagram:
\diagr{ X^\bullet=X^\bullet\otimes^{\rm ext}\ast \ar[r]\ar[d]\ar@{}[dr]|->>>>{\phantom{xx}\pushout} & X^\bullet\otimes^{\rm ext}\Delta^s/\partial\Delta^s \ar[d]& \\
        \ast \ar[r] &  X^\bullet\wedge^{\rm ext}\Delta^s/\partial\Delta^s\ar@{=}[r] & :\Sigma_{\rm ext}^sX^\bullet }
There is a dual construction called $\Omega_{\rm ext}X^\bullet$.
\end{Def}

\section{Moduli spaces in model categories}
\label{appendix:moduli}

\begin{Def} \label{haut}
Let \mc{M} be a simplicial model category and let \mc{W} be its subcategory of weak equivalences. For a cofibrant and fibrant object $X$ we define the {\bf simplicial monoid of self equivalences} denoted by {\bf haut\boldmath$\hrl X\hrr$\unboldmath} by setting
    $$ \text{haut}\hrl X\hrr_n := \Hom{\mc{W}}{X\otimes\Delta^n}{X}. $$
If we need to specify a model structure on \mc{M}, because there are several possible choices, we write the name of the structure as an index, so e.g. ${\rm haut}_{\rm Reedy}\hrl X^\bullet\hrr$ denotes the simplicial monoid in the Reedy structure on $c\mc{M}$.  
\end{Def}

\begin{Bem} \label{haut als Komponenten von map}
It is an easy, but quite important observation that the space $\text{haut}\hrl X\hrr$ consists of those connected components of the space $\map\hrl X,X\hrr$ stemming from the simplicial structure of \mc{M} that are given by the vertices corresponding to weak self equivalences.
\end{Bem}

\begin{Def} \label{Modulraum eines Objektes}
Let \mc{M} be a model category. We define the {\bf moduli space of an object} $X$ in \mc{M} to be the nerve of the following category: objects are the objects of \mc{M} that are weakly equivalent to $X$ and morphisms are the weak equivalences.
It is denoted by \boldmath$\mc{M}\hrl X\hrr$\unboldmath. Note that for each $X$ in \mc{M} this moduli space is non-empty and con\-nec\-ted.
If $S$ is a set of objects in \mc{M}, we define \boldmath$\mc{M}\hrl S\hrr$\unboldmath\ to be the nerve of the full subcategory of \mc{W}, whose objects are weakly equivalent to an element of $S$. 

We define the {\bf moduli space of a morphism} in \mc{M} in the same way as for objects: Let $f$ be an object in the category Mor$(\mc{M})$. It can be given a model structure with objectwise weak equivalences, for which we refer to \cite{DSpa:model}. 
The moduli space of $f$ is the space $\mc{M}\hrl f\hrr$ from definition \ref{Modulraum eines Objektes} in the category Mor$(\mc{M})$.
\end{Def}

The important theorem about moduli spaces is the following one proved in \cite[Prop. 2.3.]{DK:classification-for-diagrams}.
\begin{satz} \label{DK-Charakterisierung von Modulräumen}
             \label{DK-Charakterisierung von Modulraeumen}
Let \mc{M} be a simplicial model category and let $X$ be an object of \mc{M}. Then the moduli space $\mc{M}\emph{\hrl}X\emph{\hrr}$ is weakly equivalent to the space $B{\rm haut}\emph{\hrl} X\emph{\hrr}$.
\end{satz}

\section{Extension of categories}
\label{subsection:Ausdehnung}

The definitions in this paragraph are taken from \cite[VI.5]{Baues:combinatorial}. 

\begin{Def} \label{Faktorisierungskategorie}
Let \mc{C} be a category. Let $\Fac\mc{C}$ be the {\bf category of factorizations} of \mc{C}. It is the Grothendieck construction on $\mc{C}^{\rm op}\times\mc{C}$ with respect to the functor $\Hom{\mc{C}}{\frei}{\frei}$. Explicitly it has the morphisms of \mc{C} as objects, and a morphism $f\to g$ is given by a commutative diagram of the following shape:
\diagr{ \ar[d]_f & \ar[d]^g\ar[l] \\ \ar[r] & } 
\end{Def}

\begin{Def} \label{natürliche Systeme}
A {\bf natural system of abelian groups} on a category \mc{C} is a functor from $\Fac\mc{C}$ to the category ${\rm Ab}$ of abelian groups. 
\end{Def}

\begin{Bem}
There is a canonical functor $\Fac\mc{C}\to\mc{C}^{\rm op}\times\mc{C}$, sending a morphism to its source and target. 
Hence each bifunctor $\Gamma\co\mc{C}^{\rm op}\times\mc{C}\to{\rm Ab}$ induces a natural system on \mc{C}.
In this case we will also write $\Gamma\hrl X,Y\hrr$ for $\Gamma\hrl f\hrr$ if $f\co X\to Y$ is a morphism in \mc{C}.
\end{Bem}

The following definition is taken from \cite[VI(5.4)]{Baues:combinatorial}.
\begin{Def} \label{exakte Sequenz von Kategorien}
Let $\sigma:\mc{C}\to\mc{D}$ be a functor, and let $\Gamma$ and $\Xi$ be natural systems on $\mc{D}$. We write symbolically
    $$ \Gamma\toh{\gamma}\mc{C}\toh{\sigma}\mc{D}\toh{\rm ob} \Xi $$
and call this diagram an {\bf exact sequence of categories}, if the following conditions are satisfied:
\begin{punkt}
    \item
For all objects $A$ and $B$ in \mc{C} and for each morphism $f\aus\Hom{\mc{D}}{\sigma\hrl A\hrr}{\sigma\hrl B\hrr}$ there is a transitive action $\gamma$ of the group $\Gamma(f)$ on the set $\sigma^{-1}(f)\subset\Hom{\mc{C}}{A}{B}$. 
This action satisfies the {\bf linear distributivity law}:
    $$ \hrl\alpha+\wt{f}\hrr\hrl\beta+\wt{g}\hrr= f_*\alpha+g^*\beta+\wt{f}\wt{g} $$
for all $\wt{f}\aus p^{-1}(f), \wt{g}\aus p^{-1}(g), \alpha\aus G(f)$ and $\beta\aus G(g)$,, and where we have abbreviated $\gamma(\alpha,f)$ by $\alpha+f$.
    \item
For all objects $A$ and $B$ in \mc{C} and all morphisms $f:\sigma(A)\to\sigma(B)$ in \mc{D} there is an obstruction element ${\rm ob}(f)\aus \Xi(f)$ given, such that 
    $$ {\rm ob}(f)=0 $$
if and only if there exists a morphism $\wt{f}:A\to B$ with $\sigma(\wt{f})=f$.
    \item
For all $f:\sigma(A)\to\sigma(B)$ and $g:\sigma(B)\to\sigma(C)$ we have the following equation:
    $$ {\rm ob}(gf)=g_*{\rm ob}(f)+f^*{\rm ob}(g) $$
    \item
For all objects $A$ in \mc{C} and for all $\alpha\aus \Xi(\id_{\sigma(A)})$ there is an object $B$ in \mc{C} with the property that $\sigma(A)=\sigma(B)$ and ${\rm ob}(\id_{\sigma(A)})=\alpha$.
\end{punkt}
\end{Def}

The next lemma follows from the axioms and is taken from \cite[V(5.7)]{Baues:combinatorial}. 
\begin{lemma} \label{sigma entdeckt Isos}
If the functor $\sigma$ is part of an exact sequence of categories as in \emph{\ref{exakte Sequenz von Kategorien}}, then it detects isomorphism.
\end{lemma}

\end{document}